\documentclass[12pt,leqno]{article}
\usepackage{amssymb}
\usepackage[mathscr]{eucal}
\usepackage{amsmath,amssymb,latexsym,theorem,bbm}
\usepackage{amsmath,amssymb,latexsym,theorem}
\usepackage{color}
\usepackage{appendix}

\newcommand{\NN}{\mathbb{N}}

\newcommand{\RR}{\mathbb{R}}

\newcommand{\ZZ}{\mathbb{Z}}

\newcommand{\bA}{{\boldsymbol{A}}}

\newcommand{\bB}{{\boldsymbol{B}}}

\newcommand{\tb}{\widetilde{b}}

\newcommand{\tbB}{\widetilde{\bB}}

\newcommand{\bc}{{\boldsymbol{c}}}
\newcommand{\bC}{{\boldsymbol{C}}}

\newcommand{\bD}{{\boldsymbol{D}}}

\newcommand{\be}{{\boldsymbol{e}}}

\newcommand{\Bf}{{\boldsymbol{f}}}

\newcommand{\bF}{{\boldsymbol{F}}}

\newcommand{\br}{{\boldsymbol{r}}}

\newcommand{\bv}{{\boldsymbol{v}}}

\newcommand{\bx}{{\boldsymbol{x}}}
\newcommand{\bX}{{\boldsymbol{X}}}
\newcommand{\by}{{\boldsymbol{y}}}
\newcommand{\bY}{{\boldsymbol{Y}}}

\newcommand{\bz}{{\boldsymbol{z}}}

\newcommand{\bw}{{\boldsymbol{w}}}
\newcommand{\bW}{{\boldsymbol{W}}}
\newcommand{\Bbeta}{{\boldsymbol{\beta}}}
\newcommand{\tbeta}{\widetilde{\beta}}
\newcommand{\tBbeta}{\widetilde{\Bbeta}}

\newcommand{\blambda}{{\boldsymbol{\lambda}}}

\newcommand{\bxi}{{\boldsymbol{\xi}}}
\newcommand{\bmu}{{\boldsymbol{\mu}}}

\newcommand{\bbeta}{{\boldsymbol{\beta}}}

\newcommand{\bzero}{{\boldsymbol{0}}}

\newcommand{\cA}{{\mathcal A}}
\newcommand{\cB}{{\mathcal B}}

\newcommand{\cF}{{\mathcal F}}

\newcommand{\cR}{{\mathcal R}}

\newcommand{\cU}{{\mathcal U}}

\newcommand{\cc}{\mathrm{c}}
\newcommand{\dd}{\mathrm{d}}
\newcommand{\ee}{\mathrm{e}}

\newcommand{\EE}{\operatorname{\mathbb{E}}}

\newcommand{\PP}{\operatorname{\mathbb{P}}}

\newcommand{\var}{\operatorname{Var}}

\newcommand{\tN}{\widetilde{N}}

\renewcommand{\mid}{\,|\,}

\renewcommand{\leq}{\leqslant}
\renewcommand{\geq}{\geqslant}

\newcommand{\as}{\stackrel{{\mathrm{a.s.}}}{\longrightarrow}}

\newcommand{\bbone}{\mathbbm{1}}

\newcommand{\proofend}{\hfill\mbox{$\Box$}}

\numberwithin{equation}{section}

\theoremstyle{change} \theorembodyfont{\em}
\newtheorem{Lem}{Lemma.}[section]
\newtheorem{Thm}[Lem]{Theorem.}
\newtheorem{Pro}[Lem]{Proposition.}
\newtheorem{Cor}[Lem]{Corollary.}
\newtheorem{Def}[Lem]{Definition.}

\theorembodyfont{\rm}
\newtheorem{Rem}[Lem]{Remark.}

\begin{document}

\begin{center}
 {\bfseries\Large Moment formulas for multi-type continuous state and \\[2mm]
                   continuous time branching process with immigration}
 \\[6mm]

 {\sc\large
  M\'aty\'as $\text{Barczy}^{*,\diamond}$,
  \ Zenghu $\text{Li}^{**}$,
  \ Gyula $\text{Pap}^{***}$}

\end{center}

\vskip0.2cm

\noindent
 * Faculty of Informatics, University of Debrecen,
   Pf.~12, H--4010 Debrecen, Hungary.

\noindent
 ** School of Mathematical Sciences, Beijing Normal University,
     Beijing 100875,  People's Republic of China.

\noindent
 *** Bolyai Institute, University of Szeged,
     Aradi v\'ertan\'uk tere 1, H--6720 Szeged, Hungary.

\noindent e--mails: barczy.matyas@inf.unideb.hu (M. Barczy),
                    lizh@bnu.edu.cn (Z. Li),
                    papgy@math.u-szeged.hu (G. Pap).

\noindent $\diamond$ Corresponding author.

\vskip0.2cm


\renewcommand{\thefootnote}{}
\footnote{\textit{2010 Mathematics Subject Classifications\/}:
          60J80, 60J75, 60H10.}
\footnote{\textit{Key words and phrases\/}:
 multi-type continuous state and continuous time branching process with
 immigration, moments.}
\vspace*{0.2cm}
\footnote{The research of M. Barczy and G. Pap was realized in the frames of
 T\'AMOP 4.2.4.\ A/2-11-1-2012-0001 ,,National Excellence Program --
 Elaborating and operating an inland student and researcher personal support
 system''.
The project was subsidized by the European Union and co-financed by the
 European Social Fund.
Z. Li has been partially supported by NSFC under Grant No.\ 11131003
 and 973 Program under Grant No.\ 2011CB808001.}

\vspace*{-10mm}

\begin{abstract}
Recursions for moments of multi-type continuous state and continuous time
 branching process with immigration are derived.
It turns out that the \ $k$-th (mixed) moments and the \ $k$-th (mixed)
 central moments are polynomials of the initial value of the process, and
 their degree are at most \ $k$ \ and \ $\lfloor k/2 \rfloor$, \ respectively.
\end{abstract}

\section{Introduction}
\label{section_intro}

Moment formulas and estimations play an important role in the theory of
 stochastic processes, since they are useful in proving limit theorems for
 processes and for functionals of processes as well.
Branching processes form a distinguished class, since they are frequently used
 for modelling real data sets describing dynamic behaviour of populations,
 phenomenas in epidemiology, cell kinetics and genetics, so moment estimation
 for them is of great importance as well.
The main purpose of the present paper is to derive recursions for moments of a
 multi-type continuous state and continuous time branching process with
 immigration (CBI process) using the identification of such a process as a
 pathwise unique strong solution of certain stochastic differential equation
 with jumps, see \eqref{SDE_X}.

For a special Dawson--Watanabe superprocess (without immigration) with a
 special branching mechanism a recursion for the moments has been provided by
 Dynkin \cite{Dyn1} and Konno and Shiga \cite[Lemma 2.1]{KonShi}, see also
 Li \cite[Example 2.8]{Li}.
Further, Dynkin \cite[Chapter 5, Theorems 1.1 and 1.2]{Dyn2} gave
 recursive moment formulae for Dawson--Watanabe superprocesses.
We emphasize that our technique for deriving recursions for moments is
 completely different from that of Dynkin \cite{Dyn2}.
Li \cite[Propositions 2.27 and 2.38]{Li} derived formulas for the first and
 second moments for such processes.
For the class of regular immigration superprocesses, which contains multitype
 CBI processes, Li \cite[Propositions 9.11 and 9.14]{Li} derived first and
 second order moment formulas using an explicit form for the Laplace transform
 of the transition semigroup of the processes in question.

In Filipovi\'{c} et al. \cite[formula (4.4)]{FilMaySch}, one can find a formal
 representation of polynomial moments of affine processes, which include
 multitype CBI processes as well.
The idea behind this formal representation is that the infinitesimal generator
 of an affine process formally maps the finite-dimensional linear space of all
 polynomials of degree less than or equal to \ $k$ \ into itself, where
 \ $k \in \NN$, \ which suggests that the \ $k$-th moment of an affine process
 is a polynomial of the initial value of degree at most \ $k$.
\ Very recently, Filipovi\'{c} and Larsson
 \cite[Lemma 4.12 and Theorem 4.13]{FilLar} provided moment formulas and
 moment estimations for so-called polynomial preserving diffusion processes.

Yamazato \cite{Yam} considered time continuous Markov chains on the state
 space of non-negative integers having the so-called branching property and
 allowing random immigration whenever the population size is zero (as a
 special state-dependent immigration).
He investigated under which conditions the process in question has finite
 first and second moments, see \cite[Theorem 3]{Yam}, and in the so-called
 critical case he also pointed out that the first moment is a first order
 polynomial of the initial value of the process, while the second moment is a
 second order polynomial, see \cite[Theorem 5]{Yam}.

Dareiotis et al.~\cite[Lemma 2]{DarKumSab} derived some moment bounds for
 the pathwise unique strong solution of a stochastic differential equation
 (SDE) with jumps having coefficients satisfying some local Lipschitz
 condition.
We emphasize that the coefficients of the SDE of a multi-type CBI process
 given in \eqref{SDE_X} do not satisfy the locally Lipschitz condition A-5 in
 Dareiotis et al.~\cite{DarKumSab}, so their result can not be applied to a
 multi-type CBI process.
However, our technique is somewhat similar to theirs in the sense that they
 also use It\^{o}'s formula and Gronwall's inequality.

For some moment estimates for L\'evy processes, see Luschgy and Pag\`{e}s
 \cite{LusPag}; for nonlocal SDEs with time-varying delay, see Hu and Huang
 \cite{HuHua}; for linear SDEs driven by analytic fractional Brownian motion,
 see Unterberger \cite{Unt}; for unstable INteger-valued AutoRegressive models
 of order 2 (INAR(2)), see Barczy et al. \cite[Appendix A]{BarIspPap}; for a
 super-Brownian motion in one dimension with constant branching rate, see
 Perkins \cite[Lemma III.4.6]{Per}; for discrete time multi-type branching
 random walks, see  G\"un et al. \cite{GunKonSek1}, \cite{GunKonSek2}, where
 the main input comes from the many-to-few lemma due to Harris and Roberts
 \cite[Lemma 3]{HarRob}.
D\"oring and Roberts \cite[Lemma 3]{DorRob} provided a recursion for moments
 for a spatial version of a Galton--Watson process for which a system of
 branching particles moves in space and particles branch only in the presence
 of a catalyst.

The paper is organized as follows.
In Section \ref{section_CBI}, for completeness and better readability, we
 recall from Barczy et al.~\cite{BarLiPap2} some notions and statements for
 multi-type CBI processes such as the form of their infinitesimal generator,
 their branching and immigration mechanisms,
 and their representation as pathwise unique strong solutions of certain SDEs
 with jumps, see Theorem \ref{strong_solution}.
In Section \ref{section_approximation}, we consider an appropriately truncated
 version \eqref{SDE_X_K} of the SDE \eqref{SDE_X} of a multi-type CBI process,
 where we truncate the integrand of the integral with respect to a
 (non-compensated) Poisson random measure.
We show that, under some moment conditions, this truncated SDE has a pathwise
 unique strong solution which is a multi-type CBI process with explicitly
 given parameters,
 especially, the jump measures of the branching and immigration
 mechanisms are truncated, see Theorem \ref{strong_solution_K}.
Then we prove a comparison theorem with respect to the truncation mentioned
 above, see Theorem \ref{comparison_K}, and, as a consequence, we show that
 the truncated CBI process at a time point \ $t$ \ converges in \ $L^1$ \ and
 almost surely to the non-truncated CBI process at the time point \ $t$ \ as
 the level of truncation tends to \ $\infty$, \ see Theorem \ref{conv_K}.
Section \ref{section_recursion} is devoted to deriving recursion formulas for
 moments.
First, we rewrite the SDE \eqref{SDE_X} of a multi-type CBI process in a form
  which is more suitable for calculating moments.
Namely, we eliminate integrals with respect to non-compensated Poisson random
 measures, and then we perform a linear transformation in order to remove
 randomness from the drift, see Theorem \ref{SDE_transform}.
In view of Theorem \ref{conv_K}, for the proof of the recursion formula
 \eqref{SDE_truncated3} in Theorem \ref{moment_m}, it is enough to prove a
 recursion formula for moments of a truncated CBI process.
After applying It\^{o}'s formula for powers of a truncated CBI process, we
 would like to take expectations, so we have to check martingale property of
 some stochastic integrals with respect to certain compensated Poisson random
 measures.
In order to do this, by induction with respect to \ $k$, \ we prove certain
 estimates for the \ $k$-th moments of a truncated CBI process, see
 \eqref{bound_K} and \eqref{bound}.
Truncations of the jump measures of the branching and immigration
 mechanisms are needed to avoid integrability troubles when showing martingale
 property of the stopped processes \eqref{stopped_drocesses}.
It turns out that the \ $k$-th (mixed) moments and the \ $k$-th (mixed)
 central moments are polynomials of the initial value of the process, and
 their degrees are at most \ $k$ \ and \ $\lfloor k/2 \rfloor$, \ respectively,
 see Theorems \ref{moment_m} and \ref{central_moments_m}, and Corollaries
 \ref{Cor_mixed_moments} and \ref{Cor_mixed_central_moments}.
An explicit formula for the second central moment, i.e., for the variance of a
 CBI process is given in Proposition \ref{moment_formula_2}.

In a companion paper, Barczy and Pap \cite{BarPap} used the results of the
 present paper for studying the asymptotic behavior of critical irreducible
 multi-type continuous state and continuous time branching processes with
 immigration.
Further, in Barczy et al. \cite{BarKorPap} moment estimations together with
 the results in \cite{BarPap} serve as a key tool for studying asymptotic
 behavior of conditional least squares estimators of some parameters for
 \ $2$-type doubly symmetric critical irreducible CBI processes.

\section{Multi-type CBI processes}
\label{section_CBI}

Let \ $\ZZ_+$, \ $\NN$, \ $\RR$, \ $\RR_+$  \ and \ $\RR_{++}$ \ denote the set
 of non-negative integers, positive integers, real numbers, non-negative real
 numbers and positive real numbers, respectively.
For \ $x , y \in \RR$, \ we will use the notations
 \ $x \land y := \min \{x, y\} $ \ and \ $x^+:= \max \{0, x\} $.
\ By \ $\|\bx\|$ \ and \ $\|\bA\|$, \ we denote the Euclidean norm of a vector
 \ $\bx \in \RR^d$ \ and the induced matrix norm of a matrix
 \ $\bA \in \RR^{d\times d}$, \ respectively.
The natural basis in \ $\RR^d$ \ and the Borel $\sigma$-algebras on \ $\RR^d$
 \ and on \ $\RR_+^d$ \ will be denoted by \ $\be_1$, \ldots, $\be_d$, \ and
 by \ $\cB(\RR^d)$ \ and \ $\cB(\RR_+^d)$, \ respectively.
For \ $\bx = (x_i)_{i\in\{1,\ldots,d\}} \in \RR^d$ \ and
 \ $\by = (y_i)_{i\in\{1,\ldots,d\}} \in \RR^d$, \ we will use the notation
 \ $\bx \leq \by$ \ indicating that \ $x_i \leq y_i$ \ for all
 \ $i \in \{1, \ldots, d\}$.
\ By \ $C^2_\cc(\RR_+^d,\RR)$ \ we denote the set of twice continuously
 differentiable real-valued functions on \ $\RR_+^d$ \ with compact support.
Throughout this paper, we make the conventions \ $\int_a^b := \int_{(a,b]}$
 \ and \ $\int_a^\infty := \int_{(a,\infty)}$ \ for any \ $a, b \in \RR$ \ with
 \ $a < b$.

\begin{Def}\label{Def_essentially_non-negative}
A matrix \ $\bA = (a_{i,j})_{i,j\in\{1,\ldots,d\}} \in \RR^{d\times d}$ \ is called
 essentially non-negative if \ $a_{i,j} \in \RR_+$ \ whenever
 \ $i, j \in \{1,\ldots,d\}$ \ with \ $i \ne j$, \ i.e., if \ $\bA$ \ has
 non-negative off-diagonal entries.
The set of essentially non-negative \ $d \times d$ \ matrices will be denoted
 by \ $\RR^{d\times d}_{(+)}$.
\end{Def}

\begin{Def}\label{Def_admissible}
A tuple \ $(d, \bc, \Bbeta, \bB, \nu, \bmu)$ \ is called a set of admissible
 parameters if
 \renewcommand{\labelenumi}{{\rm(\roman{enumi})}}
 \begin{enumerate}
  \item
   $d \in \NN$,
  \item
   $\bc = (c_i)_{i\in\{1,\ldots,d\}} \in \RR_+^d$,
  \item
   $\Bbeta = (\beta_i)_{i\in\{1,\ldots,d\}} \in \RR_+^d$,
  \item
   $\bB = (b_{i,j})_{i,j\in\{1,\ldots,d\}} \in \RR^{d \times d}_{(+)}$,
  \item
   $\nu$ \ is a Borel measure on \ $U_d := \RR_+^d \setminus \{\bzero\}$
    \ satisfying \ $\int_{U_d} (1 \land \|\bz\|) \, \nu(\dd \bz) < \infty$,
  \item
   $\bmu = (\mu_1, \ldots, \mu_d)$, \ where, for each
    \ $i \in \{1, \ldots, d\}$, \ $\mu_i$ \ is a Borel measure on \ $U_d$
    \ satisfying
    \begin{align}\label{help_Li_page45}
      \int_{U_d}\left[ \Vert\bz\Vert \wedge \Vert\bz\Vert^2
                     + \sum_{j \in \{1, \ldots, d\} \setminus \{i\}} z_j\;
              \right] \mu_i(\dd \bz)<\infty .
    \end{align}
  \end{enumerate}
\end{Def}

\begin{Rem}
Our Definition \ref{Def_admissible} of the set of admissible parameters is a
 special case of Definition 2.6 in Duffie et al.~\cite{DufFilSch}, which is
 suitable for all affine processes.
Further, for all \ $i\in\{1,\ldots,d\}$, \ the condition \eqref{help_Li_page45}
 is equivalent to
 \begin{align}\label{help_Li_page45_equiv}
  \int_{U_d}
     \left[ (1 \land z_i)^2
             + \sum_{j \in \{1, \ldots, d\} \setminus \{i\}} (1 \land z_j) \right]
      \mu_i(\dd \bz)
    < \infty
   \quad \text{and} \quad
   \int_{U_d} \Vert\bz\Vert \bbone_{\{\Vert\bz\Vert\geq 1\}}\,\mu_i(\dd \bz)<\infty ,
 \end{align}
 see Barczy et al.~\cite[Remark 2.3]{BarLiPap2}.
\proofend
\end{Rem}

\begin{Thm}\label{CBI_exists}
Let \ $(d, \bc, \Bbeta, \bB, \nu, \bmu)$ \ be a set of admissible parameters
 in the sense of Definition \ref{Def_admissible}.
Then there exists a unique conservative transition semigroup
 \ $(P_t)_{t\in\RR_+}$ \ acting on
 the Banach space (endowed with the supremum norm) of real-valued bounded
 Borel-measurable functions on the state space \ $\RR_+^d$ \ such that its
 infinitesimal generator is
 \begin{equation}\label{CBI_inf_gen}
  \begin{aligned}
   (\cA f)(\bx)
   &= \sum_{i=1}^d c_i x_i f_{i,i}''(\bx)
      + \langle \Bbeta + \bB \bx, \Bf'(\bx) \rangle
      + \int_{U_d} \bigl( f(\bx + \bz) - f(\bx) \bigr) \, \nu(\dd \bz) \\
   &\phantom{\quad}
      + \sum_{i=1}^d
         x_i
         \int_{U_d}
          \bigl( f(\bx + \bz) - f(\bx) - f'_i(\bx) (1 \land z_i) \bigr)
          \, \mu_i(\dd \bz)
  \end{aligned}
 \end{equation}
 for \ $f \in C^2_\cc(\RR_+^d,\RR)$ \ and \ $\bx \in \RR_+^d$, \ where \ $f_i'$
\ and \ $f_{i,i}''$, \ $i \in \{1, \ldots, d\}$, \ denote the first and second
 order partial derivatives of \ $f$ \ with respect to its \ $i$-th variable,
 respectively, and \ $\Bf'(\bx) := (f_1'(\bx), \ldots, f_d'(\bx))^\top$.
\ Moreover, the Laplace transform of the transition semigroup
 \ $(P_t)_{t\in\RR_+}$ \ has a representation
 \begin{align*}
  \int_{\RR_+^d} \ee^{- \langle \blambda, \by \rangle} P_t(\bx, \dd \by)
  = \ee^{- \langle \bx, \bv(t, \blambda) \rangle - \int_0^t \psi(\bv(s, \blambda)) \, \dd s} ,
  \qquad \bx \in \RR_+^d, \quad \blambda \in \RR_+^d , \quad t \in \RR_+ ,
 \end{align*}
 where, for any \ $\blambda \in \RR_+^d$, \ the continuously differentiable function
 \ $\RR_+ \ni t \mapsto \bv(t, \blambda)
    = (v_1(t, \blambda), \ldots, v_d(t, \blambda))^\top \in \RR_+^d$
 \ is the unique locally bounded solution to the system of differential
 equations
 \begin{equation}\label{EES}
   \partial_t v_i(t, \blambda) = - \varphi_i(\bv(t, \blambda)) , \qquad
   v_i(0, \blambda) = \lambda_i , \qquad i \in \{1, \ldots, d\} ,
 \end{equation}
 with
 \[
   \varphi_i(\blambda)
   := c_i \lambda_i^2 -  \langle \bB \be_i, \blambda \rangle
      + \int_{U_d}
         \bigl( \ee^{- \langle \blambda, \bz \rangle} - 1
                + \lambda_i (1 \land z_i) \bigr)
         \, \mu_i(\dd \bz)
 \]
 for \ $\blambda \in \RR_+^d$ \ and \ $i \in \{1, \ldots, d\}$, \ and
 \[
   \psi(\blambda)
   := \langle \bbeta, \blambda \rangle
      - \int_{U_d}
         \bigl( \ee^{- \langle \blambda, \bz \rangle} - 1 \bigr)
         \, \nu(\dd \bz) , \qquad
   \blambda \in \RR_+^d .
 \]
Further, the function \ $\RR_+\times\RR_+^d\ni(t,\blambda) \mapsto \bv(t, \blambda)$
 \ is continuous.
\end{Thm}

\begin{Rem}
This theorem is a special case of Theorem 2.7 of Duffie et
 al.~\cite{DufFilSch} with \ $m = d$, \ $n = 0$ \ and zero killing rate.
\proofend
\end{Rem}

\begin{Def}\label{Def_CBI}
A conservative Markov process with state space \ $\RR_+^d$ \ and with transition
 semigroup \ $(P_t)_{t\in\RR_+}$ \ given in Theorem \ref{CBI_exists} is called a
 multi-type CBI process with parameters \ $(d, \bc, \Bbeta, \bB, \nu, \bmu)$.
\ The function
 \ $\RR_+^d \ni \blambda
    \mapsto (\varphi_1(\blambda), \ldots, \varphi_d(\blambda))^\top \in \RR^d$
 \ is called its branching mechanism, and the function
 \ $\RR_+^d \ni \blambda \mapsto \psi(\blambda) \in \RR_+$ \ is called its
 immigration mechanism.
The measures \ $\mu_i$, \ $i \in \{1, \ldots, d\}$, \ and \ $\nu$ \ are the jump
 measures of the branching and immigration mechanisms, respectively.
\end{Def}

Let \ $(\bX_t)_{t\in\RR_+}$ \ be a multi-type CBI process with parameters
 \ $(d, \bc, \Bbeta, \bB, \nu, \bmu)$ \ such that \ $\EE(\|\bX_0\|) < \infty$
 \ and the moment  condition
 \begin{equation}\label{moment_condition_1}
  \int_{U_d} \|\bz\| \bbone_{\{\|\bz\|\geq1\}} \, \nu(\dd \bz) < \infty
 \end{equation}
 holds.
Then, by Lemma 3.4 in Barczy et al. \cite{BarLiPap2},
 \begin{equation}\label{EXbX}
  \EE(\bX_t) = \ee^{t\tbB} \EE(\bX_0)
               + \int_0^t \ee^{u\tbB} \tBbeta \, \dd u ,
  \qquad t \in \RR_+ ,
 \end{equation}
 where
 \begin{gather}
  \tbB := (\tb_{i,j})_{i,j\in\{1,\ldots,d\}} , \qquad
  \tb_{i,j} := b_{i,j}
              + \int_{U_d} (z_i - \delta_{i,j})^+ \, \mu_j(\dd \bz) ,
  \label{tbB} \\
  \tBbeta := \Bbeta + \int_{U_d} \bz \, \nu(\dd \bz) ,
  \label{tBbeta}
 \end{gather}
 with \ $\delta_{i,j}:=1$ \ if \ $i=j$, \ and \ $\delta_{i,j}:=0$ \ if
 \ $i\ne j$.
\ We also introduce the modified parameters
 \ $\bD := (d_{i,j})_{i,j\in\{1,\ldots,d\}}$ \ given by
 \begin{equation}\label{dij}
  d_{i,j} := \tb_{i,j} - \int_{U_d} z_i \bbone_{\{\|\bz\|\geq1\}} \, \mu_j(\dd\bz) .
 \end{equation}
Note that \ $\tbB \in \RR^{d \times d}_{(+)}$, \ $\tBbeta \in \RR_+^d$ \ and
 \ $\bD \in \RR^{d \times d}_{(+)}$, \ since
 \begin{equation}\label{help}
  \int_{U_d} \|\bz\| \, \nu(\dd\bz) < \infty , \qquad
  \int_{U_d} (z_i - \delta_{i,j})^+ \, \mu_j(\dd \bz) < \infty , \quad
  i, j \in \{1, \ldots, d\} ,
 \end{equation}
 see Barczy et al. \cite[Section 2]{BarLiPap2}.

Let \ $\cR := \bigcup_{j=0}^d \cR_j$, \ where
 \ $\cR_j$, \ $j \in \{0, 1, \ldots, d\}$, \ are disjoint sets given by
 \[
   \cR_0 := U_d \times \{ (\bzero, 0) \}^d
         \subset \RR_+^d \times (\RR_+^d \times \RR_+)^d ,
 \]
 and
 \[
   \cR_j := \{\bzero\} \times H_{j,1} \times \cdots \times H_{j,d}
         \subset \RR_+^d \times (\RR_+^d \times \RR_+)^d , \qquad
   j \in \{1, \ldots, d\} ,
 \]
 where
 \[
   H_{j,i} := \begin{cases}
              U_d \times U_1 & \text{if \ $i = j$,} \\
              \{ (\bzero, 0) \} & \text{if \ $i \ne j$.}
             \end{cases}
 \]
(Recall that \ $U_1 = \RR_{++}$.)
Let \ $m$ \ be the uniquely defined measure on
 \ $V := \RR_+^d \times (\RR_+^d \times \RR_+)^d$ \ such that
 \ $m(V \setminus \cR) = 0$ \ and its restrictions on \ $\cR_j$,
 \ $j \in \{0, 1, \ldots, d\}$, \ are
 \begin{equation}\label{m}
  m|_{\cR_0}(\dd\br) = \nu(\dd\br) , \qquad
  m|_{\cR_j}(\dd\bz, \dd u) = \mu_j(\dd\bz) \, \dd u ,
  \quad j \in \{1, \ldots, d\} ,
 \end{equation}
 where we identify \ $\cR_0$ \ with \ $U_d$ \ and \ $\cR_1$, \ldots, $\cR_d$
 \ with \ $U_d \times U_1$ \ in a natural way.
Using again this identification, let \ $f : \RR^d \times V \to \RR_+^d$, \ and
 \ $g : \RR^d \times V \to \RR_+^d$, \ be defined by
 \[
   f(\bx, \br)
   := \begin{cases}
       \bz \bbone_{\{\|\bz\|<1\}} \bbone_{\{u \leq x_j\}} ,
        & \text{if \ $\bx = (x_1, \ldots, x_d)^\top \in \RR^d$,
                \ $\br = (\bz, u) \in \cR_j$, \ $j \in \{1, \ldots, d\}$,} \\
       \bzero , & \text{otherwise,}
      \end{cases}
 \]
 \[
   g(\bx, \br)
   := \begin{cases}
       \br, & \text{if \ $\bx \in \RR^d$, \ $\br \in \cR_0$,} \\
       \bz \bbone_{\{\|\bz\|\geq1\}} \bbone_{\{u \leq x_j\}} ,
        & \text{if \ $\bx = (x_1, \ldots, x_d)^\top \in \RR^d$,
                \ $\br = (\bz, u) \in \cR_j$, \ $j \in \{1, \ldots, d\}$,} \\
       \bzero , & \text{otherwise.}
      \end{cases}
 \]
Consider the disjoint decomposition \ $\cR = V_0 \cup V_1$, \ where
 \ $V_0 := \bigcup_{j=1}^d \cR_{j,0}$ \ and
 \ $V_1 := \cR_0 \cup \bigl( \bigcup_{j=1}^d \cR_{j,1} \bigr)$ \ are disjoint
 decompositions with
 \ $\cR_{j,k} := \{\bzero\} \times H_{j,1,k} \times \cdots \times H_{j,d,k}$,
 \ $j \in \{1, \ldots, d\}$, \ $k \in \{0, 1\}$, \ and
 \[
   H_{j,i,k} := \begin{cases}
                U_{d,k} \times U_1 & \text{if \ $i = j$,} \\
                \{ (\bzero, 0) \} & \text{if \ $i \ne j$,}
               \end{cases} \qquad
   U_{d,k} := \begin{cases}
              \{ \bz \in U_d : \|\bz\| < 1 \} & \text{if \ $k = 0$,} \\
              \{ \bz \in U_d : \|\bz\| \geq 1 \} & \text{if \ $k = 1$.}
             \end{cases}
 \]
Note that \ $f(\bx, \br) = \bzero$ \ if \ $\br \in V_1$,
 \ $g(\bx, \br) = \bzero$ \ if \ $\br \in V_0$, \ hence
 \ $\be_i^\top f(\bx, \br) g(\bx, \br) \be_j = 0$ \ for all
 \ $(\bx, \br) \in \RR^d \times V$ \ and \ $i, j \in \{1, \ldots, d\}$.

Consider the following objects:
 \begin{enumerate}
  \item[(E1)]
   a probability space \ $(\Omega, \cF, \PP)$;
  \item[(E2)]
   a \ $d$-dimensional standard Brownian motion \ $(\bW_t)_{t\in\RR_+}$;
  \item[(E3)]
   a stationary Poisson point process \ $p$ \ on \ $V$ \ with characteristic
    measure \ $m$;
  \item[(E4)]
   a random vector \ $\bxi$ \ with values in \ $\RR_+^d$, \ independent of
    \ $\bW$ \ and \ $p$.
 \end{enumerate}

\begin{Rem}\label{dRM_strong}
Note that if objects (E1)--(E4) are given, then \ $\bxi$, \ $\bW$ \ and \ $p$
 \ are automatically mutually independent according to Remark 3.4 in
 Barczy et al. \cite{BarLiPap1}.
For a short review on point measures and point processes needed for this
 paper, see, e.g., Barczy et al. \cite[Section 2]{BarLiPap1}.
\proofend
\end{Rem}

Provided that the objects (E1)--(E4) are given, let
 \ $(\cF^{\bxi,\bW\!,\,p}_t)_{t\in\RR_+}$ \ denote the augmented filtration
 generated by \ $\bxi$, \ $\bW$ \ and \ $p$, \ see
 Barczy et al.\ \cite{BarLiPap1}.

Let us consider the \ $d$-dimensional SDE
 \begin{equation}\label{SDE_X}
  \begin{aligned}
   \bX_t
   &= \bX_0 + \int_0^t (\Bbeta + \bD \bX_s) \, \dd s
      + \sum_{i=1}^d \be_i \int_0^t \sqrt{2 c_i X_{s,i}^+} \, \dd W_{s,i} \\
   &\quad
      + \int_0^t \int_{V_0} f(\bX_{s-}, \br) \, \tN(\dd s, \dd \br)
      + \int_0^t \int_{V_1} g(\bX_{s-}, \br) \, N(\dd s, \dd \br) , \qquad
   t \in \RR_+ ,
  \end{aligned}
 \end{equation}
 where \ $\bX_t = (X_{t,1}, \ldots, X_{t,d})^\top$, \ $\bD$ \ is defined in
 \eqref{dij}, \ $N(\dd s, \dd\br)$ \ is the counting measure of \ $p$ \ on
 \ $\RR_{++} \times V$, \ and
 \ $\tN(\dd s, \dd\br) := N(\dd s, \dd\br) - \dd s \, m(\dd\br)$.

\begin{Def}\label{Def_strong_solution2}
Suppose that the objects \textup{(E1)--(E4)} are given.
An \ $\RR_+^d$-valued strong solution of the SDE \eqref{SDE_X} on
 \ $(\Omega, \cF, \PP)$ \ and with respect to the standard Brownian motion
 \ $\bW$, \ the stationary Poisson point process \ $p$ \ and initial value
 \ $\bxi$, \ is an \ $\RR_+^d$-valued
 \ $(\cF^{\bxi,\bW\!,\,p}_t)_{t\in\RR_+}$-adapted c\`{a}dl\`{a}g process
 \ $(\bX_t)_{t \in \RR_+}$ \ such that \ $\PP(\bX_0 = \bxi) = 1$,
 \[
   \PP\biggl( \int_0^t
               \int_{V_0} \|f(\bX_s, \br)\|^2 \, \dd s \, m(\dd \br)
                  < \infty \biggr)
        = 1 , \qquad
   \PP\biggl( \int_0^t \int_{V_1} \|g(\bX_{s-}, \br)\| \, N(\dd s, \dd \br)
                   < \infty \biggr)
         = 1
 \]
 for all \ $t \in \RR_+$, \ and equation \eqref{SDE_X} holds \ $\PP$-a.s.
\end{Def}

Note that the integrals \ $\int_0^t (\Bbeta + \bD \bX_s) \, \dd s$ \ and
 \ $\int_0^t \sqrt{2 c_i X_{s,i}^+} \, \dd W_{s,i}$, \ $i \in \{1, \ldots, d\}$,
 \ exist, since \ $\bX$ \ is c\`{a}dl\`{a}g.
For the following result see Theorem 4.6 and Remark 3.2 in
 Barczy et al.\ \cite{BarLiPap2}.

\begin{Thm}\label{strong_solution}
Let \ $(d, \bc, \Bbeta, \bB, \nu, \bmu)$ \ be a set of admissible parameters
 such that the moment condition \eqref{moment_condition_1} holds.
Suppose that objects \textup{(E1)--(E4)} are given.
If \ $\EE(\|\bxi\|) < \infty$, \ then there is a pathwise unique
 \ $\RR_+^d$-valued strong solution to the SDE \eqref{SDE_X} with initial value
 \ $\bxi$, \ and the solution is a CBI process with parameters
 \ $(d, \bc, \Bbeta, \bB, \nu, \bmu)$.
\ Moreover, for each \ $t \in \RR_+$,
 \[
   \EE\biggl( \int_0^t \int_{V_0}
               \|f(\bX_s, \br)\|^2 \, \dd s \, m(\dd \br) \biggr)
   < \infty , \qquad
   \EE\biggl( \int_0^t \int_{V_1}
               \|g(\bX_s, \br)\| \, \dd s \, m(\dd \br) \biggr)
   < \infty .
 \]
\end{Thm}

\section{Approximation of multi-type CBI processes}
\label{section_approximation}

First we study an appropriately truncated version of the SDE \eqref{SDE_X}.

\begin{Thm}\label{strong_solution_K}
Let \ $(d, \bc, \Bbeta, \bB, \nu, \bmu)$ \ be a set of admissible parameters
 such that the moment condition \eqref{moment_condition_1} holds.
Suppose that objects \textup{(E1)--(E4)} are given.
Let \ $K \in (1, \infty]$.
\ If \ $\EE(\|\bxi\|) < \infty$, \ then there is a pathwise unique
 \ $\RR_+^d$-valued strong solution to the SDE
 \begin{equation}\label{SDE_X_K}
  \begin{aligned}
   \bX_t
   &= \bX_0 + \int_0^t \bigl( \Bbeta + \bD \bX_s \bigr) \, \dd s
      + \sum_{i=1}^d \be_i \int_0^t \sqrt{2c_i X_{s,i}^+} \, \dd W_{s,i} \\
   &\quad
      + \int_0^t \int_{V_0} f(\bX_{s-}, \br) \, \tN(\dd s, \dd \br)
      + \int_0^t \int_{V_1} g_K(\bX_{s-}, \br) \, N(\dd s, \dd \br) ,
   \qquad t \in \RR_+ ,
  \end{aligned}
 \end{equation}
 with initial value \ $\bxi$, \ where the function
 \ $g_K : \RR^d \times V \to \RR_+^d$ \ is defined by
 \[
   g_K(\bx, \br)
   := \begin{cases}
       \br \bbone_{\{\|\br\|<K\}} ,
        & \text{if \ $\bx \in \RR^d$, \ $\br \in \cR_0$,} \\
       \bz \bbone_{\{1\leq\|\bz\|<K\}} \bbone_{\{u \leq x_j\}} ,
        & \text{if \ $\bx = (x_1, \ldots, x_d)^\top \in \RR^d$,} \\
        & \text{\phantom{if \ }$\br = (\bz, u) \in \cR_{j,1}$,
                \ $j \in \{1, \ldots, d\}$,} \\
       \bzero , & \text{otherwise,}
      \end{cases}
 \]
 and the solution is a CBI process with parameters
 \ $(d, \bc, \Bbeta, \bB_K, \nu_K, \bmu_K)$, \ where
 \ $\bB_K = (b_{K,i,j})_{i,j\in\{1,\ldots,d\}}$, \ $\nu_K$ \ and
 \ $\bmu_K = (\mu_{K,1}, \ldots, \mu_{K,d})$ \ are given by
 \begin{gather}\label{bKij}
  b_{K,i,j}
  := b_{i,j}
      - \delta_{i,j}
        \int_{U_d}
         (z_i \land 1) \bbone_{\{\|\bz\|\geq K\}} \, \mu_j(\dd\bz) , \\
  \nu_K(\dd\br) := \bbone_{\{\|\br\|<K\}} \, \nu(\dd\br) , \qquad
  \mu_{K,i}(\dd\bz) := \bbone_{\{\|\bz\|<K\}} \, \mu_i(\dd\bz) . \nonumber
 \end{gather}
\end{Thm}

\noindent
\textbf{Proof.} \
In case of \ $K = \infty$, \ the SDE \eqref{SDE_X_K} coincides with the SDE
 \eqref{SDE_X}, since \ $g_\infty = g$, \ hence, by Theorem
 \ref{strong_solution}, the SDE \eqref{SDE_X_K} with \ $K = \infty$ \ admits a
 pathwise unique \ $\RR_+^d$-valued strong solution with initial value
 \ $\bxi$, \ and the solution is a CBI process with parameters
 \ $(d, \bc, \Bbeta, \bB, \nu, \bmu)$.

For each \ $K \in (1, \infty)$,
 \begin{align}\label{fg}
  \begin{split}
   \int_0^t \int_{V_0} f(\bX_{s-}, \br) \, \tN(\dd s, \dd\br)
   &= \int_0^t \int_{V_0} f(\bX_{s-}, \br) \, \tN_K(\dd s, \dd\br) , \\
   \int_0^t \int_{V_1} g_K(\bX_{s-}, \br) \, N(\dd s, \dd\br)
   &= \int_0^t \int_{V_1} g(\bX_{s-}, \br) \, N_K(\dd s, \dd\br) ,
  \end{split}
 \end{align}
 where \ $N_K(\dd s, \dd\br)$ \ is the counting measure of the stationary
 Poisson point process \ $p_K$, \ where \ $p_K$ \ denotes the thinning of \ $p$
 \ onto \ $V_0 \cup \cR_{0,K} \cup \left(\cup_{j=1}^d \cR_{j,1,K}\right)$ \ given
 by
 \begin{gather*}
  \cR_{0,K} := \{ \br \in U_d : \|\br\| < K \} \times \{ (\bzero, 0) \}^d
  \subset \RR_+^d \times (\RR_+^d \times \RR_+)^d , \\
  \cR_{j,1,K} := \{ \bzero \} \times H_{j,1,1,K} \times \cdots \times H_{j,d,1,K}
  \subset \RR_+^d \times (\RR_+^d \times \RR_+)^d ,
  \qquad j \in \{ 1, \ldots, d \} ,
 \end{gather*}
 where
 \[
   H_{j,i,1,K}
   := \begin{cases}
       \{ \bz \in U_d : 1 \leq \|\bz\| < K \} \times U_1
        & \text{if \ $i = j$,} \\
       \{ (\bzero, 0) \} & \text{if \ $i \ne j$,}
      \end{cases}
 \]
 and \ $\tN_K(\dd s, \dd\br) := N_K(\dd s, \dd\br) - \dd s \, m_K(\br)$,
 \ where \ $m_K$ \ denotes the restriction of \ $m$ \ onto
 \ $V_0 \cup \cR_{0,K} \cup \left(\cup_{j=1}^d \cR_{j,1,K}\right)
    = \cR_{0,K} \cup \bigl(\cup_{j=1}^d (\cR_{j,0} \cup \cR_{j,1,K}) \bigr)$.
\ Note that the characteristic measure of \ $p_K$ \ is \ $m_K$ \ (this can be
 checked calculating the corresponding conditional Laplace transforms, see
 Ikeda and Watanabe \cite[page 44]{IkeWat}).
Moreover, \ $m_K|_{V_0}(\dd\br) = m|_{V_0}(\dd\br)$,
 \ $m_K|_{\cR_{0,K}}(\dd\br) = \nu_K(\dd\br)$ \ and
 \begin{align*}
  &m_K|_{\cR_{j,0}\cup\cR_{j,1,K}}(\dd\bz, \dd u)
   = m_K|_{\cR_{j,0}}(\dd\bz, \dd u) + m_K|_{\cR_{j,1,K}}(\dd\bz, \dd u) \\
  &= \bbone_{\{\|\bz\|<1\}} \mu_j(\dd\bz) \, \dd u
     + \bbone_{\{1\leq\|\bz\|<K\}} \mu_j(\dd\bz) \, \dd u
   = \mu_{K,j}(\dd\bz) \, \dd u \qquad \text{for \ $j \in \{1, \ldots, d\}$.}
 \end{align*}
Consequently, the SDE \eqref{SDE_X_K} can be rewritten as
 \begin{align}\label{SDE_X_K_mod}
  \begin{split}
   \bX_t
   &= \bX_0 + \int_0^t \bigl( \Bbeta + \bD \bX_s \bigr) \, \dd s
      + \sum_{i=1}^d \be_i \int_0^t \sqrt{2c_i X_{s,i}^+} \, \dd W_{s,i} \\
   &\quad
      + \int_0^t \int_{V_0} f(\bX_{s-}, \br) \, \tN_K(\dd s, \dd\br)
      + \int_0^t \int_{V_1} g(\bX_{s-}, \br) \, N_K(\dd s, \dd\br) ,
   \qquad t \in \RR_+ .
  \end{split}
 \end{align}
Further, for each \ $K \in (1, \infty)$, \ $\nu_K$ \ and \ $\bmu_K$ \ satisfy
 parts (v) and (vi) of Definition \ref{Def_admissible}, respectively.
Further, \ $\bB_K \in \RR_{(+)}^{d\times d}$, \ hence
 \ $(d, \bc, \Bbeta, \bB_K, \nu_K, \bmu_K)$ \ is a set of admissible
 parameters.
By Theorem \ref{strong_solution}, the SDE \eqref{SDE_X_K_mod} admits a
 pathwise unique \ $\RR_+^d$-valued strong solution with initial value
 \ $\bxi$, \ and the solution is a CBI process with parameters
 \ $(d, \bc, \Bbeta, \bB_K, \nu_K, \bmu_K)$, \ since, using \eqref{tbB} and
 \eqref{dij},
 \begin{align*}
  d_{K,i,j}
  &:= b_{K,i,j} + \int_{U_d} (z_i - \delta_{i,j})^+ \, \mu_{K,j}(\dd\bz)
   - \int_{U_d} z_i \bbone_{\{\|\bz\|\geq1\}} \, \mu_{K,j}(\dd\bz) \\
  &= b_{i,j}
     - \delta_{i,j}
       \int_{U_d} (z_i \land 1) \bbone_{\{\|\bz\|\geq K\}} \, \mu_j(\dd\bz) \\
  &\quad
     + \int_{U_d} (z_i - \delta_{i,j})^+ \bbone_{\{\|\bz\|<K\}} \, \mu_j(\dd\bz)
     - \int_{U_d} z_i \bbone_{\{1\leq\|\bz\|<K\}} \, \mu_j(\dd\bz)\\
  &= b_{i,j} + \int_{U_d} (z_i - \delta_{i,j})^+ \, \mu_j(\dd\bz)
     - \int_{U_d} z_i \bbone_{\{\|\bz\|\geq1\}} \, \mu_j(\dd\bz)
     - \int_{U_d}
        (z_i - \delta_{i,j})^+ \bbone_{\{\|\bz\|\geq K\}} \, \mu_j(\dd\bz) \\
  &\quad
     + \int_{U_d} z_i \bbone_{\{\|\bz\|\geq K\}} \, \mu_j(\dd\bz)
     - \delta_{i,j}
       \int_{U_d} (z_i \land 1) \bbone_{\{\|\bz\|\geq K\}} \, \mu_j(\dd\bz)
 \end{align*}
 equals \ $d_{i,j}$ \ for all \ $i, j \in \{1, \ldots, d\}$, \ since the sum of
 the last three terms is 0.
Especially,
 \[
   \EE\biggl( \int_0^t \int_{V_0}
               \|f(\bX_s, \br)\|^2 \, \dd s \, m_K(\dd \br) \biggr)
   < \infty , \qquad
   \EE\biggl( \int_0^t \int_{V_1}
               \|g(\bX_s, \br)\| \, \dd s \, m_K(\dd \br) \biggr)
   < \infty
 \]
 for all \ $t \in \RR_+$.
\ Using \eqref{fg}, we conclude
 \[
   \EE\biggl( \int_0^t \int_{V_0}
               \|f(\bX_s, \br)\|^2 \, \dd s \, m(\dd \br) \biggr)
   < \infty , \qquad
   \EE\biggl( \int_0^t \int_{V_1}
               \|g_K(\bX_s, \br)\| \, \dd s \, m(\dd \br) \biggr)
   < \infty .
 \]
Hence the SDE \eqref{SDE_X_K} also admits a pathwise unique
 \ $\RR_+^d$-valued strong solution with initial value \ $\bxi$, \ and the
 solution is a CBI process with parameters
 \ $(d, \bc, \Bbeta, \bB_K, \nu_K, \bmu_K)$.
\proofend

Next we prove a comparison theorem for the SDE \eqref{SDE_X_K} in \ $K$.

\begin{Thm}\label{comparison_K}
Let \ $(d, \bc, \Bbeta, \bB, \nu, \bmu)$ \ be a set of admissible parameters
 such that the moment condition \eqref{moment_condition_1} holds.
Suppose that objects \textup{(E1)--(E3)} are given.
Let \ $\bxi$ \ and \ $\bxi'$ \ be random vectors with values in \ $\RR_+^d$
 \ independent of \ $\bW$ \ and \ $p$ \ such that \ $\EE(\|\bxi\|) < \infty$,
 \ $\EE(\|\bxi'\|) < \infty$ \ and \ $\PP(\bxi \leq \bxi') = 1$.
\ Let \ $K, K' \in (1, \infty]$ \ with \ $K \leq K'$.
\ Let \ $(\bX_t)_{t\in\RR_+}$ \ be a pathwise unique \ $\RR_+^d$-valued strong
 solution to the SDE \eqref{SDE_X_K} with initial value \ $\bxi$.
\ Let \ $(\bX_t')_{t\in\RR_+}$ \ be a pathwise unique \ $\RR_+^d$-valued strong
 solution to the SDE \eqref{SDE_X_K} with initial value \ $\bxi'$ \ and with
 \ $K$ \ replaced by \ $K'$.
\ Then \ $\PP(\text{$\bX_t \leq \bX_t'$ for all $t \in \RR_+$}) = 1$.
\end{Thm}

\noindent
\textbf{Proof.} \
We follow the ideas of the proof of Theorem 3.1 of Ma \cite{Ma}, which is an
 adaptation of that of Theorem 5.5 of Fu and Li \cite{FuLi}.
There is a sequence \ $\phi_k : \RR \to \RR_+$, \ $k \in \NN$, \ of twice
 continuously differentiable functions such that
 \begin{enumerate}
  \item[(i)]
   $\phi_k(z) \uparrow z^+$ \ as \ $k \to \infty$ \ for all \ $z\in\RR$;
  \item[(ii)]
   $\phi_k'(z) \in [0, 1]$ \ for all \ $z \in \RR_+$ \ and \ $k \in \NN$;
  \item[(iii)]
   $\phi_k'(z) = \phi_k(z) = 0$ \ whenever \ $-z \in \RR_+$ \ and
    \ $k \in \NN$;
  \item[(iv)]
   $\phi_k''(x - y) (\sqrt{x} - \sqrt{y})^2 \leq 2/k$ \ for all
    \ $x, y \in \RR_+$ \ and \ $k \in \NN$.
 \end{enumerate}
For a construction of such functions, see, e.g., the proof of Theorem 3.1 of
 Ma \cite{Ma}.
Let \ $\bY_t = (Y_{t,1}, \ldots, Y_{t,d})^\top := \bX_t - \bX_t'$ \ for all
 \ $t \in \RR_+$.
\ By the SDE \eqref{SDE_X_K}, we have
 \begin{align*}
  Y_{t,i}
  &= Y_{0,i}
     + \int_0^t \be_i^\top \bD \bY_s \, \dd s
     + \int_0^t
        \sqrt{2c_i} \Bigl(\sqrt{X_{s,i}} - \sqrt{X_{s,i}'}\Bigr) \, \dd W_{s,i} \\
   &\quad
      + \int_0^t \int_{V_0}
         \be_i^\top ( f(\bX_{s-}, \br) - f(\bX_{s-}', \br) )
         \, \tN(\dd s, \dd \br) \\
   &\quad
      + \int_0^t \int_{V_1}
         \be_i^\top ( g_K(\bX_{s-}, \br) - g_{K'}(\bX_{s-}', \br) )
         \, N(\dd s, \dd \br)
 \end{align*}
 for all \ $t \in \RR_+$ \ and \ $i \in \{1, \ldots, d\}$.
\ For each \ $m \in \NN$, \ put
 \[
   \tau_m := \inf\Bigl\{ t \in \RR_+
                         : \max_{i \in \{1, \ldots, d\}} \max\{X_{t,i}, X_{t,i}'\}
                           \geq m \Bigr\} .
 \]
By It\^o's formula, we obtain
 \[
   \phi_k(Y_{t\land\tau_m,i}) = \phi_k(Y_{0,i}) + \sum_{\ell=1}^7 I_{i,m,k,\ell}(t)
 \]
 for all \ $t \in \RR_+$, \ $i \in \{1, \ldots, d\}$ \ and \ $k, m \in \NN$,
 \ where
 \begin{align*}
  &I_{i,m,k,1}(t)
   := \int_0^{t\land\tau_m}
       \phi_k'(Y_{s,i}) \bigl( \be_i^\top \bD \bY_s \bigr) \, \dd s , \\
  &I_{i,m,k,2}(t)
   := \int_0^{t\land\tau_m}
       \phi_k'(Y_{s,i}) \sqrt{2c_i} \Bigl(\sqrt{X_{s,i}} - \sqrt{X_{s,i}'}\Bigr)
       \, \dd W_{s,i} , \\
  &I_{i,m,k,3}(t)
   := \frac{1}{2}
      \int_0^{t\land\tau_m}
       \phi_k''(Y_{s,i}) 2 c_i \Bigl(\sqrt{X_{s,i}} - \sqrt{X_{s,i}'}\Bigr)^2
       \, \dd s , \\
  &I_{i,m,k,4}(t)
   := \int_0^{t\land\tau_m} \int_{V_0}
       \Bigl[\phi_k\bigl(Y_{s-,i}
                         + \be_i^\top
                           ( f(\bX_{s-}, \br) - f(\bX_{s-}', \br) )\bigr)
             - \phi_k(Y_{s-,i})\Bigr]
        \tN(\dd s, \dd\br) , \\
  &I_{i,m,k,5}(t)
   := \int_0^{t\land\tau_m} \int_{V_0}
       \Bigl[ \phi_k\bigl(Y_{s-,i}
                           + \be_i^\top
                             ( f(\bX_{s-}, \br) - f(\bX_{s-}', \br) ) \bigr)
               - \phi_k(Y_{s-,i}) \\
  &\phantom{I_{i,m,k,5}(t)
            := \int_0^{t\land\tau_m} \int_{V_0} \Bigl[}
               - \phi_k'(Y_{s-,i})
                 \be_i^\top ( f(\bX_{s-}, \br) - f(\bX_{s-}', \br) ) \Bigr]
        \dd s \, m(\dd\br) , \\
  &I_{i,m,k,6}(t)
   := \int_0^{t\land\tau_m} \int_{V_1}
       \Bigl[\phi_k\bigl(Y_{s-,i}
                          + \be_i^\top
                            ( g_K(\bX_{s-}, \br) - g_{K'}(\bX_{s-}', \br) )\bigr)
              - \phi_k(Y_{s-,i})\Bigr]
       N(\dd s, \dd\br) .
 \end{align*}
Using formula (3.8) in Chapter II in Ikeda and Watanabe \cite{IkeWat}, the
 last integral can be written as
 \ $I_{i,m,k,6}(t) = I_{i,m,k,7}(t) + I_{i,m,k,8}(t)$, \ where
 \begin{align*}
  I_{i,m,k,7}(t)
  &:= \int_0^{t\land\tau_m} \int_{V_1}
       \Bigl[\phi_k\bigl(Y_{s-,i}
                          + \be_i^\top
                            ( g_K(\bX_{s-}, \br) - g_{K'}(\bX_{s-}', \br) )\bigr)
              - \phi_k(Y_{s-,i})\Bigr]
       \tN(\dd s, \dd\br) , \\
  I_{i,m,k,8}(t)
  &:= \int_0^{t\land\tau_m} \int_{V_1}
       \Bigl[\phi_k\bigl(Y_{s-,i}
                          + \be_i^\top
                            ( g_K(\bX_{s-}, \br) - g_{K'}(\bX_{s-}', \br) )\bigr)
              - \phi_k(Y_{s-,i})\Bigr]
       \dd s \, m(\dd\br) ,
 \end{align*}
 since the function
 \begin{equation}\label{I_{i,m,k,6,1}}
   \RR_+ \times V \times \Omega \ni (s, \br, \omega)
   \mapsto \phi_k\bigl(Y_{s-,i}(\omega)
                       + \be_i^\top
                         ( g_K(\bX_{s-}(\omega), \br)
                           - g_{K'}(\bX_{s-}'(\omega), \br) )\bigr)
           - \phi_k(Y_{s-,i}(\omega))
 \end{equation}
 belongs to the class \ $\bF_p^1$ \ for each \ $i \in \{1, \ldots, d\}$
 \ defined on page 62 in Ikeda and Watanabe \cite{IkeWat}.
Indeed, the predictability follows from part (iii) of Lemma A.1 in
 Barczy et al.~\cite{BarLiPap1}, and
 \begin{align*}
  &\EE\biggl( \int_0^{t\land\tau_m} \int_{V_1}
               \Bigl|\phi_k\bigl(Y_{s-,i}
                                 + \be_i^\top
                                   (g(\bX_{s-}, \br) - g_K(\bX_{s-}', \br))\bigr)
                     - \phi_k(Y_{s-,i})\Bigr|
               \dd s \, m(\dd\br) \biggr) \\
  &\qquad\qquad
   \leq \EE\biggl( \int_0^{t\land\tau_m} \int_{V_1}
                    \bigl| \be_i^\top
                           (g_K(\bX_{s-}, \br) - g_{K'}(\bX_{s-}', \br)) \bigr|
                    \, \dd s \, m(\dd\br) \biggr) ,
 \end{align*}
 where we used that by properties (ii) and (iii) of the function \ $\phi_k$,
 \ we have \ $\phi_k'(u) \in [0, 1]$ \ for all \ $u \in \RR$, \ and hence, by
 mean value theorem,
 \begin{equation}\label{(ii)}
  - z \leq \phi_k(y - z) - \phi_k(y) \leq 0 \leq \phi_k(y + z) - \phi_k(y)
      \leq z ,
  \qquad y \in \RR , \quad z \in \RR_+ , \quad k \in \NN .
 \end{equation}
We have
 \ $\be_i^\top ( g_K(\bX_{s-}, \br) - g_{K'}(\bX_{s-}', \br) )
    = r_i (\bbone_{\{\|\br\|<K\}} - \bbone_{\{\|\br\|<K'\}})
    = - r_i \bbone_{\{K\leq\|\br\|<K'\}}$
 \ for \ $\br \in \cR_0$, \ and
 \begin{equation}\label{g-gK}
  \begin{aligned}
   &\be_i^\top ( g_K(\bX_{s-}, \br) - g_{K'}(\bX_{s-}', \br) )
    = z_i (\bbone_{\{1\leq\|\bz\|<K\}} \bbone_{\{u\leq X_{s-,j}\}}
           - \bbone_{\{1\leq\|\bz\|<K'\}} \bbone_{\{u\leq X_{s-,j}'\}}) \\
   &= \begin{cases}
       z_i  & \text{if \ $Y_{s-,j} > 0$,
                    \ $X_{s-,j}' < u \leq X_{s-,j}$ \ and
                    \ $1 \leq \|\bz\| < K$,} \\
       -z_i & \text{if \ $Y_{s-,j} < 0$,
                    \ $X_{s-,j} < u \leq X_{s-,j}'$ \ and
                    \ $1 \leq \|\bz\| < K$,} \\
            & \text{or if \ $u \leq X_{s-,j}'$ \ and
                    \ $K \leq \|\bz\| < K'$,} \\
         0  & \text{otherwise,}
      \end{cases} \qquad \text{$\PP$-a.s.}
  \end{aligned}
 \end{equation}
 for \ $\br = (\bz, u) \in \cR_j$, \ $j \in \{1, \ldots, d\}$.
\ Consequently,
 \begin{align*}
  &\EE\biggl( \int_0^{t\land\tau_m} \int_{V_1}
               \Bigl|\phi_k\bigl(Y_{s-,i}
                                 + \be_i^\top
                                   (g_K(\bX_{s-}, \br)
                                    - g_{K'}(\bX_{s-}', \br))\bigr)
                     - \phi_k(Y_{s-,i})\Bigr|
               \, \dd s \, m(\dd\br) \biggr) \\
  &\leq\EE\biggl( \int_0^{t\land\tau_m} \int_{U_d}
                 r_i \bbone_{\{K\leq\|\br\|<K'\} }
                 \, \dd s \, \nu(\dd\br) \biggr) \\
  &\quad
     + \sum_{j=1}^d
        \EE\biggl( \int_0^{t\land\tau_m} \int_{U_d} \int_{U_1}
                    z_i \bbone_{\{X_{s-,j}'<u\leq X_{s-,j}\}}
                    \bbone_{\{Y_{s-,j}>0\}} \bbone_{\{1\leq\|\bz\|<K\}}
                    \, \dd s \, \mu_j(\dd\bz) \, \dd u \biggr)
   \end{align*}
 \begin{align*}
  &\quad
     + \sum_{j=1}^d
        \EE\biggl( \int_0^{t\land\tau_m} \int_{U_d} \int_{U_1}
                    z_i \bbone_{\{X_{s-,j}<u\leq X_{s-,j}'\}}
                    \bbone_{\{Y_{s-,j}<0\}} \bbone_{\{1\leq\|\bz\|<K\}}
                    \, \dd s \, \mu_j(\dd\bz) \, \dd u \biggr) \\
  &\quad
     + \sum_{j=1}^d
        \EE\biggl( \int_0^{t\land\tau_m} \int_{U_d} \int_{U_1}
                    z_i \bbone_{\{u\leq X_{s-,j}'\}} \bbone_{\{K\leq\|\bz\|<K'\}}
                    \, \dd s \, \mu_j(\dd\bz) \, \dd u \biggr) \\
  &= \EE\biggl( \int_0^{t\land\tau_m} \int_{U_d}
                    r_i \bbone_{\{K\leq\|\br\|<K'\}}
                    \, \dd s \, \nu(\dd\br) \biggr) \\
  &\quad + \sum_{j=1}^d
        \EE\biggl( \int_0^{t\land\tau_m} \int_{U_d}
                    z_i \bbone_{\{1\leq\|\bz\|<K\}} Y_{s-,j}
                    \bbone_{\{Y_{s-,j}>0\}}
                    \, \dd s \, \mu_j(\dd\bz) \biggr) \\
  &\quad + \sum_{j=1}^d
        \EE\biggl( \int_0^{t\land\tau_m} \int_{U_d}
                    z_i \bbone_{\{1\leq\|\bz\|<K\}} (-Y_{s-,j})
                    \bbone_{\{Y_{s-,j}<0\}}
                    \, \dd s \, \mu_j(\dd\bz) \biggr) \\
  &\quad + \sum_{j=1}^d
        \EE\biggl( \int_0^{t\land\tau_m} \int_{U_d}
                    z_i \bbone_{\{K\leq\|\bz\|<K'\}} X_{s-,j}'
                    \, \dd s \, \mu_j(\dd\bz) \biggr) \\
  &\leq t \int_{U_d} \|\br\| \bbone_{\{\|\br\|\geq1\}} \, \nu(\dd\br)
        + 2 m t \sum_{j=1}^d
                 \int_{U_d} \|\bz\| \bbone_{\{\|\bz\|\geq1\}} \, \mu_j(\dd\bz)
   < \infty
 \end{align*}
 by the moment condition \eqref{moment_condition_1} and
 \eqref{help_Li_page45}.

As in the proof of Lemma 4.2 in Barczy et al.\ \cite{BarLiPap2}, we obtain
 that the processes \ $(I_{i,m,k,2}(t))_{t\in\RR_+}$ \ and
 \ $(I_{i,m,k,4}(t))_{t\in\RR_+}$ \ are martingales.
Moreover, the process \ $\bigl( I_{i,m,k,7}(t)\bigr)_{t\in\RR_+}$ \ is also a
 martingale by Ikeda and Watanabe \cite[page 62]{IkeWat}, since the function
 \eqref{I_{i,m,k,6,1}} belongs to the class \ $\bF_p^1$.

Using that the matrix \ $\bD$ \ has non-negative off-diagonal entries and
 properties (ii) and (iii) of the function \ $\phi_k$, \ we obtain
 \begin{align*}
  I_{i,m,k,1}(t)
  &= \int_0^{t\land\tau_m}
      \phi_k'(Y_{s,i})
      \biggl( d_{i,i} Y_{s,i}
              + \sum_{j\in\{1,\ldots,d\}\setminus\{i\}} d_{i,j} Y_{s,j} \biggr)
      \bbone_{\RR_+}(Y_{s,i}) \, \dd s \\
  &= \int_0^{t\land\tau_m}
      \phi_k'(Y_{s,i})
      \biggl( d_{i,i} Y_{s,i}^+
              + \sum_{j\in\{1,\ldots,d\}\setminus\{i\}}
                 d_{i,j} Y_{s,j}\bbone_{\RR_+}(Y_{s,i}) \biggr)
      \dd s \\
  &\leq \int_0^{t\land\tau_m}
         \biggl( |d_{i,i}| Y_{s,i}^+
                 + \sum_{j\in\{1,\ldots,d\}\setminus\{i\}} d_{i,j} Y_{s,j}^+ \biggr)
         \, \dd s
   = \sum_{j=1}^d |d_{i,j}| \int_0^{t\land\tau_m} Y_{s,j}^+ \, \dd s .
 \end{align*}
By property (iv) of the function \ $\phi_k$,
 \[
   I_{i,m,k,3}(t) \leq (t\land\tau_m) c_i \frac{2}{k} \leq \frac{2 c_i t}{k} .
 \]
As in the proof of Lemma 4.2 in Barczy et al.\ \cite{BarLiPap2}, by
 \eqref{help}, we obtain
 \[
   I_{i,m,k,5}(t)
   \leq \frac{t}{k} \int_{U_d} z_i^2 \bbone_{\{\|\bz\|<1\}} \, \mu_i(\dd\bz)
        + \sum_{j\in\{1,\ldots,d\}\setminus\{i\}}
           \int_0^{t\land\tau_m} Y_{s,j}^+ \, \dd s
           \int_{U_d} z_i \, \mu_j(\dd\bz) .
 \]
Using again \eqref{g-gK} and integrating with respect to the variable \ $u$,
 \ we get \ $I_{i,m,k,8}(t) = \sum_{\ell=1}^4 I_{i,m,k,8,\ell}(t)$, \ where
 \begin{align*}
  &I_{i,m,k,8,1}(t)
   := \int_0^{t\land\tau_m} \! \! \! \int_{U_d}
       \bigl[ \phi_k (Y_{s-,i} - r_i) - \phi_k(Y_{s-,i}) \bigr]
       \bbone_{\{K\leq\|\br\|<K'\}}
        \, \dd s \, \nu(\dd\br) , \\
  &I_{i,m,k,8,2}(t)
   := \sum_{j=1}^d
       \int_0^{t\land\tau_m} \! \! \! \int_{U_d}
        \bigl[ \phi_k (Y_{s-,i} + z_i) - \phi_k(Y_{s-,i}) \bigr]
        Y_{s-,j} \bbone_{\{Y_{s-,j}>0\}} \bbone_{\{1\leq\|\bz\|<K\}}
        \, \dd s \, \mu_j(\dd\bz) , \\
  &I_{i,m,k,8,3}(t)
   := \sum_{j=1}^d
       \int_0^{t\land\tau_m} \! \! \! \int_{U_d}
        \bigl[ \phi_k (Y_{s-,i} - z_i) - \phi_k(Y_{s-,i}) \bigr]
        (-Y_{s-,j}) \bbone_{\{Y_{s-,j}<0\}} \bbone_{\{1\leq\|\bz\|<K\}}
        \, \dd s \, \mu_j(\dd\bz) , \\
  &I_{i,m,k,8,4}(t)
   := \sum_{j=1}^d
       \int_0^{t\land\tau_m} \! \! \! \int_{U_d}
        \bigl[ \phi_k (Y_{s-,i} - z_i) - \phi_k(Y_{s-,i}) \bigr]
        X_{s-,j}' \bbone_{\{K\leq\|\bz\|<K'\}}
        \, \dd s \, \mu_j(\dd\bz) .
 \end{align*}
By \eqref{help_Li_page45},
 \ $\int_{U_d} z_i \bbone_{\{\|\bz\|\geq1\}} \, \mu_j(\dd z) < \infty$ \ for all
 \ $i, j \in \{1, \ldots, d\}$, \ thus applying \eqref{(ii)}, we obtain
 \begin{align*}
  I_{i,m,k,8,2}(t)
  &= \sum_{j=1}^d
      \int_0^{t\land\tau_m} \int_{U_d}
       \Bigl[ \phi_k(Y_{s-,i} + z_i) - \phi_k(Y_{s-,i}) \Bigr]
       Y_{s-,j}^+ \bbone_{\{1\leq\|\bz\|<K\}} \, \dd s \, \mu_j(\dd\bz)
 \end{align*}
 \begin{align*}
  &\leq \sum_{j=1}^d
         \int_0^{t\land\tau_m} Y_{s,j}^+ \, \dd s
         \int_{U_d} z_i \bbone_{\{\|\bz\|\geq1\}} \, \mu_j(\dd\bz) .
 \end{align*}
By \eqref{(ii)}, we obtain \ $I_{i,m,k,8,1}(t) \leq 0$,
 \ $I_{i,m,k,8,3}(t) \leq 0$ \ and \ $I_{i,m,k,8,4}(t) \leq 0$.

Summarizing, we have
 \begin{align*}
  \phi_k(Y_{t\land\tau_m,i})
  &\leq \phi_k(Y_{0,i})
        + C_i \sum_{j=1}^d \int_0^{t\land\tau_m} Y_{s,j}^+ \, \dd s
        + \frac{2 c_i t}{k}
        + \frac{t}{k}
          \int_{U_d} z_i^2 \bbone_{\{\|\bz\|<1\}} \, \mu_i(\dd\bz) \\
  &\quad
        + I_{i,m,k,2}(t) + I_{i,m,k,4}(t) + I_{i,m,k,6,1}(t)
 \end{align*}
 for all \ $t \in \RR_+$, \ where
 \begin{align*}
  C_i := \max_{j\in\{1,\ldots,d\}} |d_{i,j}|
         + \max_{j\in\{1,\ldots,d\}\setminus\{i\}} \int_{U_d} z_i \, \mu_j(\dd\bz)
         + \int_{U_d} z_i \bbone_{\{\|\bz\|\geq1\}} \, \mu_i(\dd\bz).
 \end{align*}
The proof can be completed exactly as in the proof of Lemma 4.2 in
 Barczy et al.\ \cite{BarLiPap2} using Gronwall's inequality.
\proofend

Next we give a useful approximation for a multi-type CBI process.

\begin{Thm}\label{conv_K}
Let \ $(d, \bc, \Bbeta, \bB, \nu, \bmu)$ \ be a set of admissible parameters
 such that the moment condition \eqref{moment_condition_1} holds.
Suppose that objects \textup{(E1)--(E4)} are given with
 \ $\EE(\|\bxi\|) < \infty$.
\ Let \ $(\bX_t)_{t\in\RR_+}$ \ be a pathwise unique \ $\RR_+^d$-valued strong
 solution to the SDE \eqref{SDE_X} with initial value \ $\bxi$.
\ For each \ $K \in (1, \infty)$, \ let \ $(\bX_{K,t})_{t\in\RR_+}$ \ be a
 pathwise unique \ $\RR_+^d$-valued strong solution to SDE \eqref{SDE_X_K} with
 initial value \ $\bxi$.
\ Then
 \ $\PP(\text{$\bX_{K,t} \leq \bX_{K',t} \leq \bX_t$ for all $t\in\RR_+$})
    = 1$
 \ for all \ $K, K' \in (1, \infty)$ \ with \ $K \leq K'$.
\ Moreover, \ $\EE(\bX_t - \bX_{K,t}) \to \bzero$ \ and
 \ $\bX_{K,t} \uparrow \bX_t$ \ $\PP$-a.s.\ as \ $K \to \infty$ \ for all
 \ $t \in \RR_+$.
\end{Thm}

\noindent
\textbf{Proof.} \
The first statement follows from Theorem \ref{comparison_K}.
Further, by \eqref{SDE_X} and \eqref{SDE_X_K}, for each \ $K \in (1, \infty)$,
 \ $t \in \RR_+$, \ and \ $i \in \{1, \ldots, d\}$, \ we have
 \begin{equation}\label{SDE_t_Kt}
 \begin{aligned}
  X_{t,i} - X_{K,t,i}
  &= \int_0^t \be_i^\top \bD (\bX_s - \bX_{K,s}) \, \dd s
     + \int_0^t
        \sqrt{2c_i} \Bigl( \sqrt{X_{s,i}} - \sqrt{X_{K,s,i}} \Bigr)
        \dd W_{s,i} \\
  &\quad
     + \int_0^t \int_{V_0}
         \be_i^\top ( f(\bX_{s-}, \br) - f(\bX_{K,s-}, \br) )
         \, \tN(\dd s, \dd \br) \\
   &\quad
      + \int_0^t \int_{V_1}
         \be_i^\top ( g(\bX_{s-}, \br) - g_K(\bX_{K,s-}, \br) )
         \, N(\dd s, \dd \br) .
 \end{aligned}
 \end{equation}
Here
 \ $\int_0^t
     \sqrt{2c_i} \Bigl( \sqrt{X_{s,i}} - \sqrt{X_{K,s,i}} \Bigr) \dd W_{s,i}$,
 \ $t \in \RR_+$, \ is a martingale, since
 \begin{align*}
  \EE\left(\int_0^t
            2c_i
            \Bigl( \sqrt{X_{s,i}} - \sqrt{X_{K,s,i}} \Bigr)^2 \, \dd s \right)
  \leq 4c_i \int_0^t \EE(X_{s,i} + X_{K,s,i}) \, \dd s
  \leq 8c_i \int_0^t \EE(X_{s,i}) \, \dd s
  < \infty
 \end{align*}
 due to \ $\PP(\text{$\bX_{K,t} \leq \bX_t$ for all $t\in\RR_+$}) = 1$ \ and
 \eqref{EXbX}.
The process
 \[
   \int_0^t \int_{V_0}
         ( f(\bX_{s-}, \br) - f(\bX_{K,s-}, \br) )
         \, \tN(\dd s, \dd \br) ,
   \qquad t \in \RR_+ ,
 \]
 is a martingale, since the mapping
 \ $\RR_+ \times V \times \Omega \ni (s, \br, \omega)
    \mapsto f(\bX_{s-}(\omega), \br) - f(\bX_{K,s-}(\omega), \br) \in \RR^d$
 \ is in the (multidimensional versions of the) class \ $\bF_p^2$ \ defined on
 page 62 in Ikeda and Watanabe \cite{IkeWat}.
The mapping
 \ $\RR_+ \times V \times \Omega \ni (s, \br, \omega)
    \mapsto g(\bX_{s-}(\omega), \br) - g_K(\bX_{K,s-}(\omega), \br) \in \RR^d$
 \ is in the (multidimensional versions of the) class \ $\bF_p^1$, \ hence
 formula (3.8) in Chapter II in Ikeda and Watanabe \cite{IkeWat} yields
 \begin{align*}
  &\EE\biggl( \int_0^t \int_{V_1}
               \be_i^\top ( g(\bX_{s-}, \br) - g_K(\bX_{K,s-}, \br) )
               \, N(\dd s, \dd \br) \biggr) \\
  &= \EE\biggl( \int_0^t \int_{V_1}
                 \be_i^\top ( g(\bX_{s-}, \br) - g_K(\bX_{K,s-}, \br) )
                 \, \dd s \, m(\dd\br) \biggr) \\
  &= t \int_{U_d} r_i \bbone_{\{\|\br\|\geq K\}} \, \nu(\dd\br)
     + \sum_{j=1}^d
        \int_0^t \EE(X_{s,j}) \, \dd s
        \int_{U_d} z_i \bbone_{\{\|\bz\|\geq K\}} \, \mu_j(\dd\bz) \\
  &\quad
     + \sum_{j=1}^d
        \int_0^t \EE(X_{s,j} - X_{K,s,j}) \, \dd s
        \int_{U_d} z_i \bbone_{\{1\leq\|\bz\|<K\}} \, \mu_j(\dd\bz) ,
 \end{align*}
 since
 \ $\be_i^\top ( g(\bX_{s-}, \br) - g_K(\bX_{K,s-}, \br) )
    = r_i (1 - \bbone_{\{\|\br\|<K\}}) = r_i \bbone_{\{\|\br\|\geq K\}}$
 \ for \ $\br \in \cR_0$, \ and
  \begin{align*}
   &\be_i^\top ( g(\bX_{s-}, \br) - g_K(\bX_{K,s-}, \br) )
    = z_i (\bbone_{\{\|\bz\|\geq1\}} \bbone_{\{u\leq X_{s-,j}\}}
          - \bbone_{\{1\leq\|\bz\|<K\}} \bbone_{\{u\leq X_{K,s-,j}\}}) \\
   &= \begin{cases}
       z_i  & \text{if \ $X_{K,s-,j} < u \leq X_{s-,j}$
               \ and \ $1 \leq \|\bz\| < K$,} \\
            & \text{or if \ $u \leq X_{s-,j}$ \ and \ $\|\bz\| \geq K$,} \\
         0  & \text{otherwise,}
      \end{cases} \qquad \text{$\PP$-a.s.}
  \end{align*}
 for \ $\br = (\bz, u) \in \cR_j$, \ $j \in \{1, \ldots, d\}$ \ (due to
 \ $\PP(X_{K,s-,j} \leq X_{s-,j})=1$).

Hence, by taking the expectations in \eqref{SDE_t_Kt}, we obtain
 \begin{align*}
  \EE(X_{t,i} - X_{K,t,i})
  &= \int_0^t \be_i^\top \bD \EE(\bX_s - \bX_{K,s}) \, \dd s
     + \sum_{j=1}^d
        \int_0^t \EE(X_{s,j}) \, \dd s
        \int_{U_d} z_i \bbone_{\{\|\bz\|\geq K\}} \, \mu_j(\dd\bz) \\
  &\quad
     + t \int_{U_d} r_i \bbone_{\{\|\br\|\geq K\}} \, \nu(\dd\br)
     + \sum_{j=1}^d
        \int_0^t \EE(X_{s,j} - X_{K,s,j}) \, \dd s
        \int_{U_d} \!\!z_i \bbone_{\{1\leq\|\bz\|<K\}} \, \mu_j(\dd\bz) .
 \end{align*}
Thus
 \[
   \sum_{i=1}^d \EE(X_{t,i} - X_{K,t,i})
   \leq \alpha_K(t)
        + C \int_0^t \Biggl( \sum_{j=1}^d \EE(X_{s,j} - X_{K,s,j}) \Biggr) \dd s ,
 \]
 where
 \begin{align*}
  \alpha_K(t)
  &:= t \sum_{i=1}^d \int_{U_d} r_i \bbone_{\{\|\br\|\geq K\}} \, \nu(\dd\br)
      + \sum_{i=1}^d \sum_{j=1}^d
         \int_0^t \EE(X_{s,j}) \, \dd s \int_{U_d}
          z_i \bbone_{\{\|\bz\|\geq K\}} \, \mu_j(\dd\bz) ,
 \end{align*}
 \begin{align*}
  C &:= \max_{j\in\{1,\ldots,d\}}
         \sum_{i=1}^d
          \biggl( |d_{i,j}|
                  + \int_{U_d}
                     z_i \bbone_{\{\|\bz\|\geq1\}} \, \mu_j(\dd\bz) \biggr) .
 \end{align*}
By Gronwall's inequality and using that \ $\alpha_K(t)$, \ $t \in \RR_+$, \ is
 monotone increasing for each \ $K \in (1, \infty)$, \ we get
 \begin{align*}
  0 \leq \sum_{i=1}^d \EE(X_{t,i} - X_{K,t,i})
    \leq \alpha_K(t) + C\int_0^t \alpha_K(s) \ee^{C(t-s)} \, \dd s
    \leq \alpha_K(t) + \alpha_K(t)C\int_0^t \ee^{C(t-s)} \, \dd s ,
 \end{align*}
 hence \ $\EE(\bX_t - \bX_{K,t}) \to \bzero$ \ as \ $K \to \infty$ \ for all
 \ $t \in \RR_+$ \ follows from \ $\alpha_K(t) \to 0$ \ as \ $K \to \infty$
 \ (which holds by dominated convergence theorem).
Finally, a non-increasing sequence of random variables converging to 0 in
 \ $L_1$ \ automatically converges to 0 almost surely, hence
 \ $\bX_{K,t} \uparrow \bX_t$ \ $\PP$-a.s.\ as \ $K \to \infty$ \ for all
 \ $t \in \RR_+$.
\proofend

\section{Recursions for moments of multi-type CBI processes}
\label{section_recursion}

First we rewrite the SDE \eqref{SDE_X} in a form which does not contain
 integrals with respect to non-compensated Poisson random measures, and then
 we perform a linear transformation in order to remove randomness from the
 drift.
This form will be very useful in calculating moments.

\begin{Lem}\label{SDE_transform}
Let \ $(d, \bc, \Bbeta, \bB, \nu, \bmu)$ \ be a set of admissible parameters
 such that the moment condition \eqref{moment_condition_1} holds.
Suppose that objects \textup{(E1)--(E4)} are given with
 \ $\EE(\|\bxi\|) < \infty$.
\ Let \ $(\bX_t)_{t\in\RR_+}$ \ be a pathwise unique \ $\RR_+^d$-valued strong
 solution to the SDE \eqref{SDE_X} with initial value \ $\bxi$.
\ Then
 \begin{align}\label{Ito_kerdeses}
  \begin{split}
   \ee^{-t\tbB} \bX_t
   &= \bX_0 + \int_0^t \ee^{-s\tbB} \tBbeta \, \dd s
      + \sum_{k=1}^d
         \int_0^t \ee^{-s\tbB} \be_k \sqrt{2 c_k X_{s,k}} \, \dd W_{s,k} \\
   &\quad
      + \int_0^t \int_V \ee^{-s\tbB} h(\bX_{s-}, \br) \, \tN(\dd s, \dd\br) ,
   \qquad t \in \RR_+ ,
  \end{split}
 \end{align}
 where the function \ $h : \RR^d \times V \to \RR^d$ \ is defined by
 \ $h := f + g$.
\end{Lem}

\noindent
\textbf{Proof.} \
The SDE \eqref{SDE_X} can be written in the form
 \begin{align}\label{Ito_kerdeses2}
 \begin{split}
  X_{t,i}
  &= X_{0,i} + \int_0^t \be_i^\top \bigl( \Bbeta + \bD \bX_s \bigr) \, \dd s
     + \int_0^t \sqrt{2c_i X_{s,i}} \, \dd W_{s,i}
     + \int_0^t \int_{\cR_0} r_i \, N(\dd s, \dd\br) \\
  &\quad
     + \sum_{j=1}^d
        \int_0^t \int_{\cR_{j,0}}
         z_i \bbone_{\{u\leq X_{s-,j}\}} \, \tN(\dd s, \dd\br)
     + \sum_{j=1}^d
        \int_0^t \int_{\cR_{j,1}}
         z_i \bbone_{\{u\leq X_{s-,j}\}} \, N(\dd s, \dd\br)
  \end{split}
 \end{align}
 for \ $t \in \RR_+$ \ and \ $i \in \{1, \ldots, d\}$.
\ Using formula (3.8) in Chapter II in Ikeda and Watanabe \cite{IkeWat}, for
 each \ $j \in \{1, \ldots, d\}$,
 \begin{align*}
  \int_0^t \int_{\cR_{j,1}} z_i \bbone_{\{u\leq X_{s-,j}\}} \, N(\dd s, \dd\br)
  &= \int_0^t \int_{\cR_{j,1}}
      z_i \bbone_{\{u\leq X_{s-,j}\}} \, \tN(\dd s, \dd\br) \\
  &\quad
     + \int_0^t \int_{U_d} \int_{U_1}
        z_i \bbone_{\{\|\bz\|\geq1\}} \bbone_{\{u\leq X_{s-,j}\}}
        \, \dd s \, \mu_j(\dd\bz) \, \dd u,
 \end{align*}
 since
 \[
   \int_0^t \int_{U_d} \int_{U_1}
    z_i \bbone_{\{\|\bz\|\geq1\}} \bbone_{\{u\leq X_{s-,j}\}}
    \, \dd s \, \mu_j(\dd\bz) \, \dd u
   = \int_0^t X_{s,j} \, \dd s
     \int_{U_d} z_i \bbone_{\{\|\bz\|\geq1\}} \, \mu_j(\dd\bz) ,
 \]
 and consequently
 \[
   \EE\left( \int_0^t X_{s,j} \, \dd s
             \int_{U_d} z_i \bbone_{\{\|\bz\|\geq1\}} \, \mu_j(\dd\bz) \right)
   = \int_0^t \EE(X_{s,j}) \, \dd s
     \int_{U_d} z_i \bbone_{\{\|\bz\|\geq1\}} \, \mu_j(\dd\bz)
   < \infty .
 \]
In a similar way,
 \[
   \int_0^t \int_{\cR_0} r_i \, N(\dd s, \dd\br)
   = \int_0^t \int_{\cR_0} r_i \, \tN(\dd s, \dd\br)
     + \int_0^t \int_{U_d} r_i \, \dd s \, \nu(\dd\br) ,
 \]
 since
 \[
   \int_0^t \int_{U_d} r_i \, \dd s \, \nu(\dd\br)
   = t \int_{U_d} r_i \, \nu(\dd\br)
   < \infty .
 \]
Consequently, by \eqref{tBbeta},
 \begin{align}\label{help_referee}
  \begin{aligned}
   \bX_t
   &= \bX_0 + \int_0^t \bigl( \tBbeta + \tbB \bX_s \bigr) \, \dd s
      + \sum_{i=1}^d \be_i^\top \int_0^t \sqrt{2c_i X_{s,i}} \, \dd W_{s,i} \\
   &\quad
      + \int_0^t \int_V h(\bX_{s-}, \br) \, \tN(\dd s, \dd\br)
  \end{aligned}
 \end{align}
 for \ $t \in \RR_+$, \ since, by \eqref{dij},
 \begin{align*}
  &\be_i^\top \bD \bX_s
   + \sum_{j=1}^d
      X_{s,j} \int_{U_d} z_i \bbone_{\{\|\bz\|\geq1\}} \, \mu_j(\dd\bz) \\
  &= \sum_{j=1}^d
      \left( d_{i,j}
             + \int_{U_d} z_i \bbone_{\{\|\bz\|\geq1\}} \, \mu_j(\dd\bz) \right)
      X_{s,j}
   = \sum_{j=1}^d \tb_{i,j} X_{s,j}
   = \be_i^\top \tbB \bX_s .
 \end{align*}
The statement of the lemma follows by an application of the multidimensional
 It\^o's formula
 (see, e.g., Ikeda and Watanabe \cite[Chapter II, Theorem 5.1]{IkeWat}).
Indeed, for each \ $i \in \{1, \ldots, d\}$,
 \ $\be_i^\top \ee^{-t\tbB} \bX_t = F_i(t, \bX_t)$ \ with the function
 \ $F_i(t, \bx) := \be_i^\top \ee^{-t\tbB} \bx
    = \sum_{j=1}^d \be_i^\top \ee^{-t\tbB} \be_j x_j$
 \ for \ $t \in \RR_+$ \ and \ $\bx = (x_1, \ldots, x_d)^\top \in \RR^d$.
\ We have
 \ $\partial_t F_i(t, \bx) = \be_i^\top \ee^{-t\tbB} (-\tbB) \bx$,
 \ $\partial_{x_k} F_i(t, \bx) = \be_i^\top \ee^{-t\tbB} \be_k$,
 \ $\partial_{x_k} \partial_{x_\ell} F_i(t, \bx) = 0$,
 \ $i, k, \ell \in \{1, \ldots, d\}$, \ hence
 \begin{align*}
  &\be_i^\top \ee^{-t\tbB} \bX_t
   = \be_i^\top \bX_0 + \int_0^t \be_i^\top \ee^{-s\tbB} (-\tbB) \bX_s \, \dd s \\
  &\hspace*{11.3mm}
   + \sum_{k=1}^d
      \int_0^t \be_i^\top \ee^{-s\tbB} \be_k \sqrt{2 c_k X_{s,k}} \, \dd W_{s,k}
   + \sum_{k=1}^d
      \int_0^t
       \be_i^\top \ee^{-s\tbB} \be_k \be_k^\top (\tBbeta + \tbB \bX_s) \, \dd s \\
  &\hspace*{11.3mm}
   + \int_0^t \int_V
      \Bigl[ \be_i^\top \ee^{-s\tbB} \bigl(\bX_{s-} + h(\bX_{s-}, \br)\bigr)
             - \be_i^\top \ee^{-s\tbB} \bX_{s-} \Bigr]
      \tN(\dd s, \dd\br) \\
  &\hspace*{11.3mm}
   + \int_0^t \int_V
      \biggl[ \be_i^\top \ee^{-s\tbB} \bigl(\bX_s + h(\bX_s, \br)\bigr)
              - \be_i^\top \ee^{-s\tbB} \bX_s
              - \sum_{k=1}^d
                 (\be_i^\top \ee^{-s\tbB} \be_k) \be_k^\top h(\bX_s, \br) \biggr]
      \dd s \, m(\dd\br) ,
 \end{align*}
 which yields the statement of the lemma (indeed, the integrand and hence the
 integral with respect to the measure \ $\dd s \, m(\dd\br)$ \ is identically
 zero).
\proofend

\begin{Rem}
We point out that in the proof of Lemma \ref{SDE_transform} formally we have
 no right to apply Theorem 5.1 in Ikeda and Watanabe
 \cite[Chapter II]{IkeWat} for \eqref{help_referee}, since the integrand
 of the integral \ $\int_0^t\int_V h(\bX_{s-},\br)\, \tN(\dd s,\dd \br)$ \ does
 not belong to the (multidimensional version of the) space \ $\bF_p^{2,loc}$.
\ Instead, we should apply It\^o's formula to \eqref{Ito_kerdeses2} (or
 equivalently to \eqref{SDE_X}).
However, after applying It\^o's formula to \eqref{Ito_kerdeses2}, one
 could rewrite the obtained equation yielding \eqref{Ito_kerdeses} under the
 moment condition \eqref{moment_condition_1}, as desired.
We will use this observation in other proofs as well later on.
\proofend
\end{Rem}

\begin{Thm}\label{moment_m}
Let \ $(\bX_t)_{t\in\RR_+}$ \ be a CBI process with parameters
 \ $(d, \bc, \Bbeta, \bB, \nu, \bmu)$ \ such that
 \ $\EE(\|\bX_0\|^q) < \infty$,
 \begin{equation}\label{moment_condition_m}
  \int_{U_d} \|\bz\|^q \bbone_{\{\|\bz\|\geq1\}} \, \nu(\dd \bz) < \infty ,
  \qquad
  \int_{U_d} \|\bz\|^q \bbone_{\{\|\bz\|\geq1\}} \, \mu_i(\dd \bz) < \infty ,
  \quad
  i \in \{1, \ldots, d\}
 \end{equation}
 with some \ $q \in \NN$.
\ Then \ $\EE(\|\bX_t\|^q) < \infty$ \ for all \ $t \in \RR_+$, \ and we have
 the recursion
 \begin{equation}\label{SDE_truncated3}
  \begin{aligned}
   \EE\left(X_{t,j}^k\right)
   &= \EE\left[(\be_j^\top \ee^{t\tbB} \bX_0)^k\right]
      + k \int_0^t
           (\be_j^\top \ee^{(t-s)\tbB} \tBbeta)
           \EE\left[(\be_j^\top \ee^{(t-s)\tbB} \bX_s)^{k-1}\right]
           \dd s \\
   &\quad
    + k (k - 1)
        \sum_{i=1}^d
         c_i
         \int_0^t
          (\be_j^\top \ee^{(t-s)\tbB} \be_i)^2
          \EE\left[(\be_j^\top \ee^{(t-s)\tbB} \bX_s)^{k-2} X_{s,i}\right]
          \dd s \\
   &\quad
    + \sum_{\ell=0}^{k-2}
       \binom{k}{\ell}
        \sum_{i=1}^d
         \int_0^t \int_{U_d}
          (\be_j^\top \ee^{(t-s)\tbB} \bz)^{k-\ell}
          \EE\bigl[(\be_j^\top \ee^{(t-s)\tbB} \bX_s)^\ell X_{s,i}\bigr]
          \dd s \, \mu_i(\dd\bz) \\
   &\quad
    + \sum_{\ell=0}^{k-2}
       \binom{k}{\ell}
        \int_0^t \int_{U_d}
         (\be_j^\top \ee^{(t-s)\tbB} \bz)^{k-\ell}
         \EE\bigl[(\be_j^\top \ee^{(t-s)\tbB} \bX_s)^\ell\bigr]
         \dd s \, \nu(\dd\bz)
  \end{aligned}
 \end{equation}
 for all \ $k \in \{1, \ldots, q\}$, \ $j \in \{1, \ldots, d\}$ \ and
 \ $t \in \RR_+$.
\ Moreover, for each \ $t \in \RR_+$, \ $k \in \{1, \ldots, q\}$ \ and
 \ $j \in \{1, \ldots, d\}$, \ there exists a polynomial
 \ $Q_{t,k,j} : \RR^d \to \RR$ \ having degree at most \ $k$ \ such that
 \begin{align}\label{help9_polinomQ}
  \EE(X_{t,j}^k) = \EE[Q_{t,k,j}(\bX_0)] , \qquad t \in \RR_+ .
 \end{align}
The coefficients of the polynomial \ $Q_{t,k,j}$ \ depend on
 \ $d$, $\bc$, $\Bbeta$, $\bB$, $\nu$, $\mu_1$, \ldots, $\mu_d$.
\end{Thm}

Note that formula \eqref{SDE_truncated3} with \ $k = 1$ \ gives back formula
 \eqref{EXbX}.

\noindent
\textbf{Proof of Theorem \ref{moment_m}.} \
In the Introduction we gave a brief sketch of the present proof.
Consider objects \textup{(E1)--(E4)} with initial value
 \ $\bxi = \by = (y_1, \ldots, y_d)^\top \in \RR_+^d$.
\ Let \ $(\bY_t)_{t\in\RR_+}$ \ be a pathwise unique \ $\RR_+^d$-valued strong
 solution to the SDE \eqref{SDE_X} with initial value \ $\by$.
\ By Theorem \ref{strong_solution}, \ $\bY$ \ is a CBI process with parameters
 \ $(d, \bc, \Bbeta, \bB, \nu, \bmu)$ \ having c\`adl\`ag trajectories.
Then the finite dimensional distributions of \ $\bX$ \ conditioned that
 \ $\bX_0 = \by$ \ and \ $\bY$ \ coincide.
Let \ $K \in (1, \infty)$, \ and let \ $(\bY_{K,t})_{t\in\RR_+}$ \ be a pathwise
 unique \ $\RR_+^d$-valued strong solution to SDE \eqref{SDE_X_K}
 (or, equivalently, to SDE \eqref{SDE_X_K_mod}) with initial value \ $\by$.
\ By Theorem \ref{strong_solution_K}, \ $(\bY_{K,t})_{t\in\RR_+}$ \ is a CBI
 process with parameters \ $(d, \bc, \Bbeta, \bB_K, \nu_K, \bmu_K)$.
\ Truncation of measures \ $\nu$ \ and \ $\mu_i$, \ $i \in \{1, \ldots, d\}$,
 \ will be needed to avoid integrability troubles when showing martingale
 property of the stopped processes \eqref{stopped_drocesses}.

The aim of the following consideration is to show by induction with respect to
 \ $k$ \ that for each \ $k \in \ZZ_+$ \ and \ $K \in (1, \infty)$ \ there
 exists a continuous function \ $f_{K,k,\by} : \RR_+ \to \RR_+$ \ such that
 \begin{equation}\label{bound_K}
  \EE(\|\bY_{K,t}\|^k) \leq f_{K,k,\by}(t), \qquad t \in \RR_+ ,
 \end{equation}
 and for each \ $k \in \{0, 1, \ldots, q\}$, \ there exists a continuous
 function \ $f_{k,\by }: \RR_+ \to \RR_+$ \ such that
 \begin{equation}\label{bound}
  \sup_{K\in(1,\infty)} \EE(\|\bY_{K,t}\|^k) \leq f_{k,\by}(t), \qquad t \in \RR_+ .
 \end{equation}
For \ $k = 0$, \ \eqref{bound_K} and \eqref{bound} are trivial.
By Lemma \ref{SDE_transform},
 \begin{equation}\label{Ito}
   \bw^\top \ee^{-t\tbB_K} \bY_{K,t}
   = \bw^\top \by
      + \int_0^t \bw^\top \ee^{-s\tbB_K} \tBbeta_K \, \dd s
      + I_{K,\bw,1}(t) +
   J_{K,\bw,1,0}(t) + J_{K,\bw,1,1}(t)
 \end{equation}
 for all \ $t \in \RR_+$, \ $\bw \in \RR^d$ \ and \ $K \in (1, \infty)$,
 \ where
 \begin{align*}
  I_{K,\bw,1}(t)
  &:= \sum_{i=1}^d
       \int_0^t
        (\bw^\top \ee^{-s\tbB_K} \be_i) \sqrt{2 c_i Y_{K,s,i}} \, \dd W_{s,i} , \\
  J_{K,\bw,1,i}(t)
  &
  := \int_0^t \int_{V_i}
       \bw^\top \ee^{-s\tbB_K} h(\bY_{K,s-}, \br) \, \tN_K(\dd s, \dd\br) , \qquad
  i \in \{0, 1\} ,
 \end{align*}
 with \ $\tN_K$ \ defined in the proof of Theorem \ref{strong_solution_K},
 \ $\tBbeta_K = (\tbeta_{K,i})_{i\in\{1,\ldots,d\}}$ \ and
 \ $\tbB_K = (\tb_{K,i,j})_{i,j\in\{1,\ldots,d\}}$ \ are given by
 \begin{align*}
  \tbeta_{K,i}
  &:= \beta_i + \int_{U_d} r_i \, \nu_K(\dd\br)
    = \beta_i + \int_{U_d} r_i \bbone_{\{\|\br\|<K\}} \, \nu(\dd\br)
    = \tbeta_i - \int_{U_d} r_i \bbone_{\{\|\br\|\geq K\}} \, \nu(\dd\br) ,
 \end{align*}
 and
 \begin{align*}
  \tb_{K,i,j}
  &:= b_{K,i,j} + \int_{U_d} (z_i - \delta_{i,j})^+ \, \mu_{K,j}(\dd\bz) \\
  &\:= b_{i,j}
       - \delta_{i,j}
         \int_{U_d} (z_i \land 1) \bbone_{\{\|\bz\|\geq K\}} \, \mu_j(\dd\bz)
       + \int_{U_d} (z_i - \delta_{i,j})^+ \bbone_{\{\|\bz\|<K\}} \, \mu_j(\dd\bz) \\
  &\:= b_{i,j} + \int_{U_d} (z_i - \delta_{i,j})^+ \, \mu_j(\dd\bz)
       - \int_{U_d} z_i \bbone_{\{\|\bz\|\geq K\}} \, \mu_j(\dd\bz)
     = \tb_{i,j} - \int_{U_d} z_i \bbone_{\{\|\bz\|\geq K\}} \, \mu_j(\dd\bz) ,
 \end{align*}
 with \ $b_{K,i,j}$ \ defined in \eqref{bKij}, where we applied the identity
 \ $(z_i \land 1) + (z_i - 1)^+ = z_i$ \ for \ $z_i \in \RR_+$.
\ By It\^o's formula, we obtain
 \begin{align}\label{SDE_truncated}
 \begin{split}
  (\bw^\top \ee^{-t\tbB_K} \bY_{K,t})^k
  &= (\bw^\top \by)^k + I_{K,\bw,k}(t)
  + J_{K,\bw,k,0}(t) + J_{K,\bw,k,1}(t) \\
  &\quad
   + k \int_0^t
        (\bw^\top \ee^{-s\tbB_K} \bY_{K,s})^{k-1} (\bw^\top \ee^{-s\tbB_K} \tBbeta_K)
        \, \dd s \\
  &\quad
   + \frac{1}{2} k (k - 1)
      \int_0^t
       (\bw^\top \ee^{-s\tbB_K} \bY_{K,s})^{k-2}
       \sum_{i=1}^d (\bw^\top \ee^{-s\tbB_K} \be_i)^2 2 c_i Y_{K,s,i}
       \, \dd s \\
  &\quad
   + \int_0^t \int_V
       \biggl[ \Bigl(\bw^\top \ee^{-s\tbB_K}
                     \bigl(\bY_{K,s} + h(\bY_{K,s}, \br)\bigr)\Bigr)^k
               - (\bw^\top  \ee^{-s\tbB_K} \bY_{K,s})^k \\
  &\phantom{\hspace*{21.1mm}}
               - k (\bw^\top \ee^{-s\tbB_K} \bY_{K,s})^{k-1}
                 \bigl(\bw^\top \ee^{-s\tbB_K} h(\bY_{K,s}, \br)\bigr)\biggr]
       \dd s \, m_K(\dd\br)
  \end{split}
 \end{align}
 for all \ $k \in \NN$ \ with \ $k \geq 2$, \ $t \in \RR_+$, \ $\bw \in \RR^d$
 \ and \ $K \in (1, \infty)$, \ where
 \begin{align*}
  I_{K,\bw,k}(t)
  &:= k \sum_{i=1}^d \int_0^t
         (\bw^\top \ee^{-s\tbB_K} \bY_{K,s})^{k-1}
         (\bw^\top \ee^{-s\tbB_K} \be_i) \sqrt{2 c_i Y_{K,s,i}}
          \, \dd W_{s,i} ,\\
  J_{K,\bw,k,i}(t)
  &:= \int_0^t \int_{V_i}
       \biggl[ \Bigl(\bw^\top \ee^{-s\tbB_K}
                    \bigl(\bY_{K,s-} + h(\bY_{K,s-}, \br)\bigr)\Bigr)^k
              - (\bw^\top \ee^{-s\tbB_K} \bY_{K,s-})^k \biggr]
       \tN_K(\dd s, \dd\br) \\
  &\: = \sum_{\ell=0}^{k-1}
        \binom{k}{\ell}
        \int_0^t \int_{V_i}
         (\bw^\top \ee^{-s\tbB_K} \bY_{K,s-})^\ell
         \bigl(\bw^\top \ee^{-s\tbB_K} h(\bY_{K,s-}, \br)\bigr)^{k-\ell}
         \, \tN_K(\dd s, \dd\br)
 \end{align*}
 for \ $i \in \{0, 1\}$.
\ For each \ $n \in \NN$, \ consider the stopping time
 \ $\tau_{K,n} := \inf\{ t \in \RR_+ : \|\bY_{K,t}\| \geq n \}$.
\ Clearly, \ $\tau_{K,n} \as \infty$ \ as \ $n \to \infty$, \ since
 \ $(\bY_{K,t})_{t\in\RR_+}$ \ has c\`adl\`ag trajectories.
The stopped processes
 \begin{align}\label{stopped_drocesses}
  \left(I_{K,\bw,k}(t\land\tau_{K,n})\right)_{t\in\RR_+}
    \qquad\text{and}\qquad
  \left(J_{K,\bw,k,i}(t\land\tau_{K,n})\right)_{t\in\RR_+} , \quad i \in \{0, 1\} ,
 \end{align}
 are martingales for all \ $k, n \in \NN$, \ $K \in (1, \infty)$ \ and
 \ $\bw \in \RR^d$.
\ Indeed,
 \begin{align*}
  &\EE\left( \int_0^{t\land\tau_{K,n}}
              (\bw^\top \ee^{-s\tbB_K} \bY_{K,s})^{2k-2}
              (\bw^\top \ee^{-s\tbB_K} \be_i)^2 Y_{K,s,i}
              \, \dd s \right)
   \leq n^{2k-1} \|\bw\|^{2k} t c(t)^{2k}
   < \infty ,
 \end{align*}
 since for all \ $t \in \RR_+$ \ and \ $s \in [0, t]$, \ we have
 \[
   \|\ee^{-s\tbB_K}\| \leq \ee^{s\|\tbB_K\|}
   \leq \exp\biggl\{t \sup_{K\in(1,\infty)} \|\tbB_K\|\biggr\}
   =: c(t)
   < \infty ,
 \]
 because, for all \ $i, j \in \{1, \ldots, d\}$, \ by monotone convergence
 theorem,
 \begin{align*}
   \tb_{K,i,j}
   = \tb_{i,j} - \int_{U_d} z_i \bbone_{\{\|\bz\|\geq K\}} \, \mu_j(\dd\bz)
   \uparrow \tb_{i,j} \qquad \text{as \ $K \to \infty$.}
 \end{align*}
Moreover, for each \ $\ell \in \{0, 1, \ldots, k-1\}$,
 \begin{align*}
  &\EE\biggl( \int_0^{t\land\tau_{K,n}} \!\!\! \int_{V_0}
               \Bigl| (\bw^\top \ee^{-s\tbB_K} \bY_{K,s-})^\ell
                      \bigl(\bw^\top \ee^{-s\tbB_K}
                            h(\bY_{K,s-}, \br)\bigr)^{k-\ell} \Bigr|^2
               \dd s \, m_K(\dd\br) \biggr) \\
  &\leq \|\bw\|^{2k} c(t)^{2k}
        \sum_{j=1}^d
         \EE\biggl( \int_0^{t\land\tau_{K,n}} \!\!\! \int_{U_d} \int_{U_1} \!
                     \|\bY_{K,s-}\|^{2\ell} \|\bz\|^{2(k-\ell)} \bbone_{\{\|\bz\|<1\}}
                     \bbone_{\{u\leq Y_{K,s-,j}\}}
                     \, \dd s \, \mu_{K,j}(\dd \bz) \, \dd u \biggr) \\
  &\leq \|\bw\|^{2k} t c(t)^{2k} n^{2\ell+1}
        \sum_{j=1}^d
         \int_{U_d} \|\bz\|^{2(k-\ell)} \bbone_{\{\|\bz\|<1\}} \, \mu_j(\dd \bz)
   < \infty
 \end{align*}
 and
 \begin{align*}
  &\EE\biggl( \int_0^{t\land\tau_{K,n}} \!\!\! \int_{V_1}
               \Bigl| (\bw^\top \ee^{-s\tbB_K} \bY_{K,s-})^\ell
                      \bigl(\bw^\top \ee^{-s\tbB_K}
                            h(\bY_{K,s-}, \br)\bigr)^{k-\ell} \Bigr|
               \, \dd s \, m_K(\dd\br) \biggr) \\
  &\leq \|\bw\|^k c(t)^k
        \EE\biggl( \int_0^{t\land\tau_{K,n}} \!\!\! \int_{U_d} \!
                    \|\bY_{K,s-}\|^\ell \|\br\|^{k-\ell}
                    \, \dd s \, \nu_K(\dd \br) \biggr) \\
  &\quad
        + \|\bw\|^k c(t)^k
          \sum_{j=1}^d
           \EE\biggl( \int_0^{t\land\tau_{K,n}} \!\!\! \int_{U_d} \int_{U_1} \!
                       \|\bY_{K,s-}\|^\ell \|\bz\|^{k-\ell} \bbone_{\{\|\bz\|\geq1\}}
                       \bbone_{\{u\leq Y_{K,s-,j}\}}
                       \, \dd s \, \mu_{K,j}(\dd \bz) \, \dd u \biggr) \\
  &\leq \|\bw\|^k t c(t)^k n^\ell
        \biggl( \int_{U_d} \|\br\|^{k-\ell} \bbone_{\{\|\br\|<K\}} \, \nu(\dd \br)
                + n \sum_{j=1}^d
                     \int_{U_d}
                      \|\bz\|^{k-\ell} \bbone_{\{1\leq\|\bz\|<K\}}
                      \, \mu_j(\dd \bz) \biggr)
   < \infty ,
 \end{align*}
 hence, by Ikeda and Watanabe
 \cite[Chapter II, Proposition 2.2 and page 62]{IkeWat}, the processes in \eqref{stopped_drocesses}
 are martingales for all \ $k, n \in \NN$,
 \ $K \in (1, \infty)$ \ and \ $\bw \in \RR^d$.
\ Here we used that for all \ $k \in \NN$ \ and
 \ $\ell \in \{0, 1, \ldots, k-1\}$,
 \begin{equation}\label{nu}
  \int_{U_d} \|\br\|^{k-\ell} \bbone_{\{\|\br\|<K\}} \, \nu(\dd\br)
  \leq \int_{U_d} \|\br\| \bbone_{\{\|\br\|<1\}} \, \nu(\dd\br)
       + K^{k-\ell} \int_{U_d} \bbone_{\{1\leq\|\br\|<K\}} \, \nu(\dd\br)
  < \infty
 \end{equation}
 due to part (v) of Definition \ref{Def_admissible},
 \begin{equation}\label{mu_i}
  \begin{aligned}
   \int_{U_d} \|\bz\|^{2(k-\ell)} \bbone_{\{\|\bz\|<1\}} \, \mu_i(\dd\bz)
   &\leq \int_{U_d} \|\bz\|^2 \bbone_{\{\|\bz\|<1\}} \, \mu_i(\dd\bz) \\
   &\leq \int_{U_d}
          \biggl( z_i^2 + \sum_{j\in\{1,\ldots,d\}\setminus\{i\}} z_j \biggr)
          \bbone_{\{\|\bz\|<1\}} \, \mu_i(\dd\bz)
    < \infty
  \end{aligned}
 \end{equation}
 due to part (vi) of Definition \ref{Def_admissible}, and
 \begin{equation}\label{mu_iK}
  \begin{aligned}
   \int_{U_d} \|\bz\|^{k-\ell} \bbone_{\{1\leq\|\bz\|<K\}} \, \mu_i(\dd\bz)
   &\leq K^{k-\ell} \int_{U_d} \bbone_{\{1\leq\|\bz\|<K\}} \, \mu_i(\dd\bz)\\
   &\leq K^{k-\ell} \int_{U_d} \|\bz\|^q \bbone_{\{\|\bz\|\geq 1\}} \, \mu_i(\dd\bz)
   < \infty
  \end{aligned}
 \end{equation}
 due to assumption \eqref{moment_condition_m}.

By replacing \ $t$ \ by \ $t\land\tau_{K,n}$ \ in \eqref{Ito} and \eqref{SDE_truncated},
 and then taking expectations on both sides of these equations, we conclude
 \[
   \EE\left[\bw^\top \ee^{-(t\land\tau_{K,n})\tbB_K} \bY_{K,t\land\tau_{K,n}}\right]
   = \bw^\top \by
     + \EE\left( \int_0^{t\land\tau_{K,n}}
                   \bw^\top \ee^{-s\tbB_K} \tBbeta_K
                   \, \dd s \right)
 \]
 and
 \begin{multline*}
  \EE\left[(\bw^\top \ee^{-(t\land\tau_{K,n})\tbB_K} \bY_{K,t\land\tau_{K,n}})^k\right] \\
  \begin{aligned}
  &= (\bw^\top \by)^k
    + k \EE\left( \int_0^{t\land\tau_{K,n}}
                   (\bw^\top \ee^{-s\tbB_K} \tBbeta_K)
                   (\bw^\top \ee^{-s\tbB_K} \bY_{K,s})^{k-1}
                   \, \dd s \right) \\
  &\quad
    + k (k - 1)
      \sum_{i=1}^d
       c_i \EE\left( \int_0^{t\land\tau_{K,n}}
                      (\bw^\top \ee^{-s\tbB_K} \be_i)^2
                      (\bw^\top \ee^{-s\tbB_K} \bY_{K,s})^{k-2}
                      \, Y_{K,s,i}
                      \, \dd s \right) \\
  &\quad
    + \sum_{\ell=0}^{k-2}
       \binom{k}{\ell}
       \EE\left( \int_0^{t\land\tau_{K,n}} \!\!\! \int_V
                  (\bw^\top \ee^{-s\tbB_K} \bY_{K,s})^\ell
                  \bigl(\bw^\top \ee^{-s\tbB_K} h(\bY_{K,s}, \br)\bigr)^{k-\ell}
                  \, \dd s \, m_K(\dd\br) \right)
  \end{aligned}
 \end{multline*}
 for all \ $k, n \in \NN$ \ with \ $k \geq 2$, \ $K \in (1, \infty)$ \ and
 \ $\bw \in \RR^d$.
\ By Fatou's lemma,
 \begin{equation}\label{Fatou}
  \begin{aligned}
   &\EE\left[(\bw^\top \ee^{-t\tbB_K} \bY_{K,t})^k\right]
    = \EE\left[\lim_{n\to\infty}
                (\bw^\top \ee^{-(t\land\tau_{K,n})\tbB_K}
                 \bY_{K,t\land\tau_{K,n}})^k\right]\\
   &\qquad
    \leq \liminf_{n\to\infty}
          \EE\left[(\bw^\top \ee^{-(t\land\tau_{K,n})\tbB_K}
                   \bY_{K,t\land\tau_{K,n}})^k\right]
    \leq \|\bw\|^k \bigl(\|\by\|^k + g_{K,k,\by}(t)\bigr)
  \end{aligned}
 \end{equation}
 with
 \begin{align*}
  g_{K,k,\by}(t)
  &:= k \|\tBbeta\| c(t)^k \int_0^t \EE(\|\bY_{K,s}\|^{k-1}) \, \dd s
      + k (k - 1) c(t)^k
        \sum_{i=1}^d c_i \int_0^t \EE(\|\bY_{K,s}\|^{k-1}) \, \dd s \\
  &\phantom{:=\:}
      + c(t)^k
        \sum_{\ell=0}^{k-2}
         \binom{k}{\ell}
         \Biggl[ \int_0^t \EE(\|\bY_{K,s}\|^\ell) \, \dd s
                 \int_{U_d} \|\bz\|^{k-\ell} \, \nu_K(\dd \bz) \\
  &\phantom{:= + c(t)^k \sum_{\ell=0}^{k-2} \binom{k}{\ell} \Biggl[}
                 + \sum_{j=1}^d
                    \int_0^t \EE(\|\bY_{K,s}\|^{\ell+1}) \, \dd s
                    \int_{U_d} \|\bz\|^{k-\ell} \, \mu_{K,j}(\dd \bz) \Biggr] .
 \end{align*}
Here we used that \ $\bzero \leq \tBbeta_K \leq \tBbeta$ \ for all
 \ $K \in (1, \infty)$,
 \[
   h(\bx, \br)
   := \begin{cases}
       \br, & \text{if \ $\bx \in \RR_+^d$, \ $\br \in \cR_0$,} \\
       \bz \bbone_{\{u \leq x_j\}} ,
        & \text{if \ $\bx = (x_1, \ldots, x_d)^\top \in \RR_+^d$,
                \ $\br = (\bz, u) \in \cR_j$, \ $j \in \{1, \ldots, d\}$,} \\
       \bzero , & \text{otherwise,}
      \end{cases}
 \]
 and hence
 \begin{align*}
  &\EE\left( \int_0^{t\land\tau_{K,n}} \!\!\! \int_V
              \|\bY_{K,s}\|^\ell \|h(\bY_{K,s}, \br)\|^{k-\ell}
              \, \dd s \, m_K(\dd\br) \right) \\
  &= \EE\left( \int_0^{t\land\tau_{K,n}} \!\!\! \int_{U_d}
                \|\bY_{K,s}\|^\ell \|\br\|^{k-\ell}
                \, \dd s \, \nu_K(\dd\br) \right) \\
  &\quad
     + \sum_{j=1}^d
        \EE\left( \int_0^{t\land\tau_{K,n}} \!\!\! \int_{U_d} \int_{U_1}
                   \|\bY_{K,s}\|^\ell \|\bz\|^{k-\ell} \bbone_{\{u \leq Y_{K,s,j}\}}
                   \, \dd s \, \mu_{K,j}(\dd\bz) \, \dd u \right)
 \end{align*}
 \begin{align*}
  &\leq \int_0^t \EE(\|\bY_{K,s}\|^\ell) \, \dd s
        \int_{U_d} \|\br\|^{k-\ell} \, \nu_K(\dd \br)
        + \sum_{j=1}^d
            \int_0^t \EE(\|\bY_{K,s}\|^{\ell+1}) \, \dd s
            \int_{U_d} \|\bz\|^{k-\ell} \, \mu_{K,j}(\dd \bz) .
 \end{align*}
If we suppose that \eqref{bound_K} holds for \ $0, 1, \ldots, k-1$ \ with
 \ $k \in \NN$ \ and for some \ $K \in (1, \infty)$, \ then \ $g_{K,k,\by}$ \ is
 a continuous function on \ $\RR_+$.
\ Note that, for each \ $\bx = (x_1, \ldots, x_d) \in \RR_+^d$ \ and \ $k \in \NN$,
 \ we have
 \begin{equation}\label{pmi}
  \|\bx\|^k \leq d^{k/2} \max_{i\in\{1,\ldots,d\}} x_i^k .
 \end{equation}
For \ $k \geq 2$, \ this is a consequence of the power mean inequality,
 for \ $k = 1$, \ this is trivial.
Choosing \ $\bw := \be_i$, \ $i \in \{1, \ldots, d\}$, \ by \eqref{pmi} and \eqref{Fatou},
 we have
 \[
   \EE(\|\ee^{-t\tbB_K} \bY_{K,t} \|^k) \leq d^{k/2} (\|\by\|^k + g_{K,k,\by}(t)) ,
   \qquad t \in \RR_+ , \quad k \in \NN , \quad K \in (1, \infty) .
 \]
Consequently,
 \begin{align*}
   \EE\left[ (Y_{K,t,i})^k \right]
   &= \EE\left[ (\be_i^\top \bY_{K,t})^k \right]
    = \EE\left[ (\be_i^\top \ee^{t\tbB_K} \ee^{-t\tbB_K} \bY_{K,t})^k \right]\\
   &\leq d^{k/2} \|\be_i^\top \ee^{t\tbB_K}\|^k (\|\by\|^k + g_{K,k,\by}(t))
    \leq d^{k/2} c(t)^k (\|\by\|^k + g_{K,k,\by}(t))
 \end{align*}
 for each \ $i \in \{1, \ldots, d\}$, \ and whence, again by \eqref{pmi},
 \begin{align*}
  \EE(\|\bY_{K,t}\|^k)
  \leq d^{k/2} \max_{i\in\{1,\ldots,d\}} \EE\left[ (Y_{K,t,i})^k \right]
  \leq d^k c(t)^k (\|\by\|^k + g_{K,k,\by}(t))
   =: f_{K,k,\by}(t) ,
 \end{align*}
 where \ $f_{K,k,\by} : \RR_+ \to \RR_+$ \ is a continuous function, hence we
 obtain \eqref{bound_K} for \ $k$ \ and \ $K$.

If we suppose that \eqref{bound} holds for \ $0,1, \ldots, k-1$ \ with
 \ $k \in \{1, \ldots, q\}$, \ then the continuity of the function \ $c$ \ and
 condition \eqref{moment_condition_m} imply the existence of a continuous
 function \ $g_{k,\by} : \RR_+ \to \RR_+$ \ such that
 \begin{equation}\label{g_m}
   \sup_{K\in(1,\infty)} g_{K,k,\by}(t) \leq g_{k,\by}(t) , \qquad t \in \RR_+ .
 \end{equation}
Namely, one can choose
 \begin{align*}
  &g_{k,\by}(t)
   := k \|\tBbeta\| c(t)^k \int_0^t f_{k-1,\by}(s) \, \dd s
      + k (k - 1) c(t)^k
        \sum_{i=1}^d c_i \int_0^t f_{k-1,\by}(s) \, \dd s \\
  &\hspace*{7mm}
      + c(t)^k
        \sum_{\ell=0}^{k-2}
         \binom{k}{\ell}
         \Biggl[ \sum_{j=1}^d
                  \int_0^t f_{\ell+1,\by}(s) \, \dd s
                  \int_{U_d} \|\bz\|^{k-\ell} \, \mu_j(\dd \bz)
      + \int_0^t f_{\ell,\by}(s) \, \dd s
        \int_{U_d} \|\br\|^{k-\ell} \, \nu(\dd \br) \Biggr],
 \end{align*}
 for \ $t \in \RR_+$, \ and the continuity of \ $g_{k,\by}$ \ is obvious, since
 \begin{equation}\label{ints}
  \begin{aligned}
   \sup_{K\in(1,\infty)} \int_{U_d} \|\bz\|^{k-\ell} \, \mu_{K,j}(\dd \bz)
   &= \int_{U_d} \|\bz\|^{k-\ell} \, \mu_j(\dd \bz) ,
   \qquad j \in \{1, \ldots, d\} , \\
   \sup_{K\in(1,\infty)} \int_{U_d} \|\br\|^{k-\ell} \, \nu_K(\dd \br)
   &= \int_{U_d} \|\br\|^{k-\ell} \, \nu(\dd \br) .
  \end{aligned}
 \end{equation}
We have \ $g_{k,\by}(t) < \infty$, \ since for all \ $k \in \{1, \ldots, q\}$
 \ and \ $\ell \in \{0, 1, \ldots, k-2\}$,
 \begin{equation}\label{nu_infty}
   \int_{U_d} \|\br\|^{k-\ell} \, \nu(\dd\br)
   \leq \int_{U_d} \|\br\| \bbone_{\{\|\br\|<1\}} \, \nu(\dd\br)
        + \int_{U_d} \|\br\|^q \bbone_{\{\|\br\|\geq1\}} \, \nu(\dd\br)
   < \infty
 \end{equation}
 due to part (v) of Definition \ref{Def_admissible} and assumption
 \eqref{moment_condition_m},
 \ $\int_{U_d} \|\bz\|^{k-\ell} \bbone_{\{\|\bz\|<1\}} \, \mu_i(\dd\bz) < \infty$
 \ can be derived as in \eqref{mu_i}, and
 \begin{equation}\label{mu_infty}
   \int_{U_d} \|\bz\|^{k-\ell} \bbone_{\{\|\bz\|\geq1\}} \, \mu_i(\dd\bz)
   \leq \int_{U_d} \|\bz\|^q \bbone_{\{\|\bz\|\geq 1\}} \, \mu_i(\dd\bz)
   < \infty
 \end{equation}
 due to assumption \eqref{moment_condition_m}.
Thus \eqref{bound} holds for \ $k$ \ with the continuous function
 \ $f_{k,\by}(t) := d^k c(t)^k \bigl(\|\by\|^k + g_{k,\by}(t)\bigr)$,
 \ $t \in \RR_+$.
\ Note that \ $f_{k,\by}(t)$ \ and \ $g_{k,\by}(t)$ \ are polynomials of
 \ $\|\by\|$ \ having degree \ $k$ \ and \ $k-1$, \ respectively.

Here we point out that \eqref{bound} may not hold for any \ $k \in \NN$,
 \ but only for \ $k \in \{0, 1, \ldots, q\}$.
\ Indeed, the integrals in \eqref{ints} are not necessarily finite, thus our
 constructions for \ $f_{k,\by}$ \ and \ $g_{k,\by}$ \ do not necessarily work.

By Theorem \ref{conv_K}, \ $\bY_{K,t} \uparrow \bY_t$ \ a.s.~as
 \ $K \to \infty$.
\ Hence \ $Y_{K,t,j}^k \uparrow Y_{t,j}^k$ \ a.s.~as \ $K \to \infty$ \ for all
 \ $j \in \{1, \ldots, d\}$, \ $k \in \NN$ \ and \ $t \in \RR_+$, \ which
 yields \ $\lim_{K\to\infty} \EE(Y_{K,t,j}^k) = \EE(Y_{t,j}^k) \in [0, \infty]$ \ by
 monotone convergence theorem.
Using \eqref{bound} with \ $k = q$, \ we obtain
 \ $\EE(Y_{t,j}^q) \in [0, \infty)$, \ $t \in \RR_+$,
 \ $j \in \{1, \ldots, d\}$, \ implying
 \ $\EE(\|\bY_t\|^q) \leq f_{q,\by}(t) < \infty$ \ for all \ $t \in \RR_+$.
\ By the tower rule for conditional expectations
 (i.e., the law of iterated expectations), it suffices to show
 \begin{equation}\label{Xq}
  \EE(\|\bX_t\|^q \mid \bX_0)
  \leq f_{q,\bX_0}(t) \quad \text{$\PP$-a.s.,}
  \qquad t \in \RR_+ ,
 \end{equation}
 since \ $f_{q,\bX_0}(t)$ \ is a polynomial of \ $\|\bX_0\|$ \ having degree
 \ $q$, \ where the conditional expectation
 \ $\EE(\|\bX_t\|^q \mid \bX_0) \in [0, \infty]$ \ is meant in the generalized
 sense, see, e.g., Stroock \cite[Theorem 5.1.6]{Str}.
In order to show \eqref{Xq}, let \ $\phi_n : \RR_+ \to \RR_+$, \ $n \in \NN$,
 \ be simple functions such that \ $\phi_n(y) \uparrow y$ \ as
 \ $n \to \infty$ \ for all \ $y \in \RR_+$.
\ Then, by the monotone convergence theorem for (generalized) conditional
 expectations, see, e.g., Stroock \cite[Theorem 5.1.6]{Str}, we obtain
 \ $\EE(\phi_n(\|\bX_t\|^q) \mid \bX_0) \uparrow \EE(\|\bX_t\|^q \mid \bX_0)$
 \ as \ $n \to \infty$ \ $\PP$-almost surely.
For each \ $B \in \cB(\RR^d)$, \ we have
 \[
   \EE(\bbone_B(\bX_t) \mid \bX_0)
   = \PP(\bX_t \in B \mid \bX_0)
   = \int_{\RR_+^d} \bbone_B(\by) \, P_t(\bX_0, \dd \by)
   \quad \text{$\PP$-a.s.,}
 \]
 hence
 \ $\EE(\phi_n(\|\bX_t\|^q) \mid \bX_0)
    = \int_{\RR_+^d} \phi_n(\|\by\|^q) \, P_t(\bX_0, \dd \by)$ \ $\PP$-almost
 surely.
By the monotone convergence theorem,
 \ $\int_{\RR_+^d} \phi_n(\|\by\|^q) \, P_t(\bX_0, \dd \by)
    \uparrow \int_{\RR_+^d} \|\by\|^q \, P_t(\bX_0, \dd \by)$
 \ as \ $n \to \infty$.
\ By \ $\EE(\|\bY_t\|^q) \leq f_{q,\by}(t) < \infty$, \ we get
 \[
   \EE(\|\bX_t\|^q \mid \bX_0)
   = \int_{\RR_+^d} \|\by\|^q \, P_t(\bX_0, \dd \by)
   \leq f_{q,\bX_0}(t) \quad \text{$\PP$-a.s.,}
 \]
 hence we conclude \eqref{Xq}.

The aim of the following discussion is to show that the processes
 \[
   \left(I_{K,\bw,k}(t)\right)_{t\in\RR_+}
   \quad \text{and}\quad
   \left(J_{K,\bw,k,i}(t)\right)_{t\in\RR_+} , \quad
    i \in \{0, 1\},
 \]
 are martingales for all \ $K \in (1, \infty)$, \ $\bw \in \RR^d$ \ and
 \ $k \in \NN$.
\ These follow similarly to the earlier discussion, since the estimates
 \eqref{bound_K} yield
 \begin{align*}
  &\EE\left( \int_0^t
              (\bw^\top \ee^{-s\tbB_K} \bY_{K,s})^{2k-2}
              (\bw^\top \ee^{-s\tbB_K} \be_i)^2
              Y_{K,s,i} \, \dd s \right)
   \leq \|\bw\|^{2k} c(t)^{2k} \int_0^t f_{K,2k-1,\by}(s) \, \dd s
   < \infty , \\
  &\EE\left( \int_0^t \int_{U_d}
              \bigl| (\bw^\top \ee^{-s\tbB_K } \bY_{K,s-})^\ell
                     (\bw^\top \ee^{-s\tbB_K } \br)^{k-\ell} \bigr|
              \, \dd s \, \nu_K(\dd\br) \right) \\
  &\qquad\qquad\qquad
   \leq \|\bw\|^k c(t)^k \int_0^t f_{K,\ell,\by}(s) \, \dd s
        \int_{U_d} \|\br\|^{k-\ell} \bbone_{\{\|\br\|<K\}} \, \nu(\dd\br)
   < \infty, \\
  &\EE\left( \int_0^t \int_{U_d} \int_{U_1}
              \bigl| (\bw^\top \ee^{-s\tbB_K } \bY_{K,s-})^\ell
                     (\bw^\top \ee^{-s\tbB_K } \bz
                     \bbone_{\{u\leq Y_{K,s-,j}\}})^{k-\ell} \bigr|^2
              \bbone_{\{\|\bz\| <1\}}
              \, \dd s \, \mu_{K,j}(\dd\bz) \, \dd u \right) \\
  &\qquad\qquad\qquad
   \leq \|\bw\|^{2k} c(t)^{2k} \int_0^t f_{K,2\ell+1,\by}(s) \, \dd s
        \int_{U_d} \|\bz\|^{2(k-\ell)} \bbone_{\{\|\bz\|<1\}} \, \mu_j(\dd\bz)
   < \infty, \\
  &\EE\left( \int_0^t \int_{U_d} \int_{U_1}
              \bigl| (\bw^\top \ee^{-s\tbB_K} \bY_{K,s-})^\ell
                     (\bw^\top \ee^{-s\tbB_K} \bz
                     \bbone_{\{u\leq Y_{K,s-,j}\}})^{k-\ell} \bigr|
              \bbone_{\{\|\bz\|\geq1\}}
              \, \dd s \, \mu_{K,j}(\dd\bz) \, \dd u \right) \\
  &\qquad\qquad\qquad
   \leq \|\bw\|^k c(t)^k \int_0^t f_{K,\ell+1,\by}(s) \, \dd s
        \int_{U_d} \|\bz\|^{k-\ell} \bbone_{\{1\leq\|\bz\|<K\}} \, \mu_j(\dd\bz)
   < \infty
 \end{align*}
 for all \ $\ell \in \{0, 1, \ldots, k-1\}$, \ where we used \eqref{nu},
 \eqref{mu_i} and \eqref{mu_iK}.
Thus, taking again expectations of \eqref{SDE_truncated} and putting
 \ $\bw = \ee^{t\tbB_K^\top} \be_j$, \ $j \in \{1, \ldots, d\}$, \ we conclude
 \begin{align}\label{SDE_truncated2}
 \begin{split}
  &\EE\left(Y_{K,t,j}^k\right)
  = (\be_j^\top \ee^{t\tbB_K} \by)^k
    + k \int_0^t
         (\be_j^\top \ee^{(t-s)\tbB_K} \tBbeta_K)
         \EE\bigl[(\be_j^\top \ee^{(t-s)\tbB_K} \bY_{K,s})^{k-1}\bigr]
         \, \dd s \\
   &\hspace*{15.3mm}
   + k (k - 1)
      \sum_{i=1}^d
       c_i
       \int_0^t
        (\be_j^\top \ee^{(t-s)\tbB_K} \be_i)^2
        \EE\bigl[(\be_j^\top \ee^{(t-s)\tbB_K} \bY_{K,s})^{k-2} Y_{K,s,i}\bigr]
        \, \dd s \\
   &\hspace*{15.3mm}
    + \sum_{\ell=0}^{k-2}
       \binom{k}{\ell}
       \sum_{i=1}^d
        \int_0^t \int_{U_d}
         (\be_j^\top \ee^{(t-s)\tbB_K} \bz)^{k-\ell}
         \EE\bigl[(\be_j^\top \ee^{(t-s)\tbB_K} \bY_{K,s})^\ell Y_{K,s,i}\bigr]
         \, \dd s \, \mu_{K,i}(\dd\bz) \\
   &\hspace*{15.3mm}
    + \sum_{\ell=0}^{k-2}
       \binom{k}{\ell}
       \int_0^t \int_{U_d}
        (\be_j^\top \ee^{(t-s)\tbB_K} \bz)^{k-\ell}
        \EE\bigl[(\be_j^\top \ee^{(t-s)\tbB_K} \bY_{K,s})^\ell\bigr]
        \, \dd s \, \nu_K(\dd\bz)
  \end{split}
 \end{align}
 for all \ $j \in \{1, \ldots, d\}$, \ $t \in \RR_+$ \ and \ $k \in \NN$
 \ with \ $k \geq 2$.

Next we show \eqref{SDE_truncated3} with \ $\bX_0 = \by$ \ for all
 \ $k \in \{1, \ldots, q\}$, \ $j \in \{1, \ldots, d\}$ \ and \ $t \in \RR_+$.
\ By monotone convergence theorem, \ $\tBbeta_K \to \tBbeta$ \ and
 \ $\tbB_K \to \tbB$ \ as \ $K \to \infty$.
\ We will show by the dominated convergence theorem that the integrals in
 \eqref{SDE_truncated2} tends to those in \eqref{SDE_truncated3} as
 \ $K \to \infty$.
\ First, we check that the integrands converge pointwise.
For all \ $t \in \RR_+$, \ $s \in[0, t]$ \ and \ $j \in \{1, \ldots, d\}$,
 \ we have
 \[
   \EE\left[(\be_j^\top \ee^{(t-s)\tbB_K} \bY_{K,s})^\ell\right]
   \to \EE\left[(\be_j^\top \ee^{(t-s)\tbB} \bY_s)^\ell\right]
        = \EE\left[(\be_j^\top \ee^{(t-s)\tbB} \bX_s)^\ell\right]
 \]
 as \ $K \to \infty$ \ for all \ $\ell \in \{1, \ldots, k-1\}$, \ and
 \[
   \EE\left[(\be_j^\top \ee^{(t-s)\tbB_K} \bY_{K,s})^\ell Y_{K,s,i}\right]
   \to \EE\left[(\be_j^\top \ee^{(t-s)\tbB} \bY_s)^\ell Y_{s,i}\right]
   = \EE\left[(\be_j^\top \ee^{(t-s)\tbB} \bX_s)^\ell X_{s,i}\right]
 \]
 as \ $K \to \infty$ \ for all \ $\ell \in \{1, \ldots, k-2\}$.
\ Indeed, \ $\EE\bigl[(\be_j^\top \ee^{(t-s)\tbB_K} \bY_{K,s})^\ell\bigr]$ \ is a
 linear combination of \ $\EE(Y_{K,s,i_1} \cdots Y_{K,s,i_\ell})$,
 \ $i_1, \ldots, i_\ell \in \{1, \ldots, d\}$.
\ By Theorem \ref{conv_K}, \ $\bY_{K,t} \uparrow \bY_t$ \ a.s.~as
 \ $K \to \infty$, \ hence
 \ $Y_{K,s,i_1} \cdots Y_{K,s,i_\ell} \uparrow Y_{s,i_1} \cdots Y_{s,i_\ell}$
 \ a.s.~as \ $K \to \infty$, \ which yields
 \ $\lim_{K\to\infty} \EE(Y_{K,s,i_1} \cdots Y_{K,s,i_\ell})
    = \EE(Y_{s,i_1} \cdots Y_{s,i_\ell}) \in [0, \infty]$
 \ by monotone convergence theorem.
Using \ $\EE(\|\bY_s\|^q) < \infty$, \ we have
 \ $\EE(Y_{s,i_1} \cdots Y_{s,i_\ell}) < \infty$, \ and we can use again
 \ $\tbB_K \to \tbB$.
\ The expectation
 \ $\EE\bigl[(\be_j^\top \ee^{(t-a)\tbB_K} \bY_{K,s})^\ell Y_{K,s,i}\bigr]$ \ can be
 handled in the same way (we only note that
 \ $\EE(Y_{s,i_1} \cdots Y_{s,i_\ell} Y_{s,i}) < \infty$).
\ Next we check that the integrands can be bounded by integrable functions
 uniformly in \ $K \in (1, \infty)$.
\ Applying \eqref{Fatou} and \eqref{g_m} with \ $t = s$ \ and
 \ $\bw = \ee^{t\tbB_K^\top} \be_j$, \ and using that
 \ $\bzero \leq \tBbeta_K \leq \tBbeta$, \ we obtain
 \[
   \sup_{K\in(1,\infty)}
    \left| (\be_j^\top \ee^{(t-s)\tbB_K} \tBbeta_K)
           \EE\left[(\be_j^\top \ee^{(t-s)\tbB_K} \bY_{K,s})^{k-1}\right] \right|
   \leq \|\tBbeta\| c(t)^k \bigl(\|\by\|^{k-1} + g_{k-1,\by}(s)\bigr)
 \]
 for all \ $t \in \RR_+$, \ $s \in[0, t]$ \ and \ $j \in \{1, \ldots, d\}$.
\ The integrals in the first sum can be handled in a similar way.
Further,
 \[
   \sup_{K\in(1,\infty)}
    \left| (\be_j^\top \ee^{(t-s)\tbB_K} \bz)^{k-\ell}
           \EE\left[(\be_j^\top \ee^{(t-s)\tbB_K} \bY_{K,s})^\ell
                    Y_{K,s,i}\right] \right|
   \leq \|\bz\|^{k-\ell} c(t)^k \bigl(\|\by\|^{\ell+1} + g_{\ell,\by}(s)\bigr)
 \]
 for all \ $t \in \RR_+$, \ $s \in[0, t]$, \ $j \in \{1, \ldots, d\}$,
 \ $\bz \in \RR_+^d$, \ $\ell \in \{0, 1, \ldots, k-2\}$ \ and
 \ $k \in \{1, \ldots, q\}$, \ where the function
 \ $\RR_+^d \ni \bz \mapsto \|\bz\|^{k-\ell}$ \ is integrable with respect to
 the measures \ $\mu_i$, \ $i \in \{1, \ldots, d\}$, \ by \eqref{mu_i} and
 \eqref{mu_infty}.
The integrals in the third sum can be handled in a similar way using
 \eqref{nu_infty}.
Hence we can apply dominated convergence theorem to obtain
 \eqref{SDE_truncated3} with \ $\bX_0 = \by$.
\ By the law of total expectation we obtain \eqref{SDE_truncated3} whenever
 \ $\EE(\|\bX_0\|^q) < \infty$.

Now we turn to prove \eqref{help9_polinomQ}.
Again by the law of total probability, it is enough to prove
 \eqref{help9_polinomQ} for \ $\bY$.
\ Using the recursion \eqref{SDE_truncated3} for \ $\bY$, \ we obtain the
 existence of suitable polynomials \ $Q_{t,k,j}$, \ $t \in \RR_+$,
 \ $k \in \{1, \ldots, q\}$, \ $j \in \{1, \ldots, d\}$, \ by induction with
 respect to \ $k$.
\ Indeed, for \ $k = 1$, \ we have
 \ $\EE(Y_{t,j}) = \be_j^\top \ee^{s\tbB} \by
                  + \int_0^s \be_j^\top \ee^{v\tbB} \tBbeta \, \dd v$,
 \ $j \in \{1, \ldots, d\}$, \ $t \in \RR_+$.
\ Now, suppose that for some \ $k \in \NN$ \ with \ $k + 1 \leq q$, \ suitable
 polynomials \ $Q_{t,1,j}$, \ldots, $Q_{t,k,j}$ \ exist for all \ $t \in \RR_+$
 \ and \ $j \in \{1, \ldots, d\}$.
\ We apply the recursion \eqref{SDE_truncated3} for \ $k + 1$.
\ Then the function \ $\RR_+^d \ni \by \mapsto (\be_j^\top \ee^{t\tbB} \by)^{k+1}$
 is a polynomial of degree at most \ $k + 1$.
\ Moreover, for each \ $\ell \in \{0, 1, \ldots, k\}$ \ and \ $s, t \in \RR_+$
 \ with \ $s \leq t$, \ the function
 \[
   \RR_+^d \ni \by
   \mapsto \EE\bigl[ (\be_j^\top \ee^{(t-s)\tbB} \bY_s)^\ell \bigr]
 \]
 is a polynomial of degree at most \ $\ell \leq k$.
\ Further, for each \ $\ell \in \{0, 1, \ldots, k - 1\}$ \ and
 \ $s, t \in \RR_+$ \ with \ $s \leq t$, \ the function
 \[
   \RR_+^d \ni \by
   \mapsto \EE\bigl[ (\be_j^\top \ee^{(t-s)\tbB} \bY_s)^\ell Y_{s,i} \bigr]
 \]
 is a polynomial of degree at most \ $\ell + 1 \leq k$.
\ Consequently, by \eqref{SDE_truncated3},
 \ $\RR_+^d \ni \by \mapsto \EE\bigl[(Y_{t,j})^{k+1}\bigr]$ \ is a polynomial of
 degree at most \ $k + 1$, \ and we conclude the existence of suitable
 polynomials \ $Q_{t,k+1,j}$ \ for all \ $t \in \RR_+$ \ and
 \ $j \in \{1, \ldots, d\}$.
\proofend

For mixed moments, we have the following corollary.

\begin{Cor}\label{Cor_mixed_moments}
Let \ $(\bX_t)_{t\in\RR_+}$ \ be a CBI process with parameters
 \ $(d, \bc, \Bbeta, \bB, \nu, \bmu)$ \ such that
 \ $\EE(\|\bX_0\|^q) < \infty$ \ and the moment conditions
 \eqref{moment_condition_m} hold with some \ $q \in \NN$.
\ Then for all \ $t\in\RR_+$, \ $k \in \{1, \ldots, q\}$ \ and
 \ $i_1, \ldots, i_k \in \{1, \ldots, d\}$, \ there exists a polynomial
 \ $Q_{t,k,i_1,\ldots,i_k} : \RR^d \to \RR$ \ having degree at most \ $k$ \ such
 that
 \begin{align*}
  \EE( X_{t,i_1} \cdots X_{t,i_k} )
  = \EE(Q_{t,k,i_1,\ldots,i_k}(\bX_0)).
 \end{align*}
The coefficients of the polynomial \ $Q_{t,k,i_1,\ldots,i_k}$ \ depend on
 \ $d$, $\bc$, $\Bbeta$, $\bB$, $\nu$, $\mu_1$, \ldots, $\mu_d$.
\end{Cor}

\noindent{\bf Proof.}
By the method of the proof of Theorem \ref{moment_m} (formally replacing
 \ $\be_j$ \ by \ $\bw \in \RR^d$ \ in \eqref{SDE_truncated3}), one can derive
 \begin{align*}
   \EE\left[\langle \bw, \bX_t \rangle^k\right]
   &= \EE\left[(\bw^\top\ee^{t\tbB} \bX_0)^k\right]
      + k \int_0^t
           (\bw^\top \ee^{(t-s)\tbB} \tBbeta)
            \EE[(\bw^\top\ee^{(t-s)\tbB} \bX_s)^{k-1}]\,\dd s \\
   &\quad + k (k - 1)
       \sum_{i=1}^d
         c_i\int_0^t
         (\bw^\top \ee^{(t-s)\tbB} \be_i)^2
         \EE\bigl[(\bw^\top \ee^{(t-s)\tbB} \bX_s )^{k-2}
                 X_{s,i}\bigr]\dd s \\
   &\quad
    + \sum_{\ell=0}^{k-2}
       \binom{k}{\ell}
        \int_0^t \int_{\cU_d}
          (\bw^\top \ee^{(t-s)\tbB} \br)^{k-\ell}
          \EE\Bigl[\bigl(\bw^\top \ee^{(t-s)\tbB} \bX_s\bigr)^\ell\Bigr]
          \dd s \, \nu(\dd\br)
 \end{align*}
 \begin{align*}
   &\quad
    + \sum_{\ell=0}^{k-2}
       \binom{k}{\ell}
         \sum_{i=1}^d
          \int_0^t \int_{\cU_d}
           (\bw^\top \ee^{(t-s)\tbB} \bz)^{k-\ell}
           \EE\Bigl[\bigl(\bw^\top \ee^{(t-s)\tbB} \bX_s\bigr)^\ell
                    X_{s,i}\Bigr]
           \dd s \, \mu_i(\dd\bz)
 \end{align*}
 for all \ $t \in \RR_+$, \ $k \in \{1, \ldots, q\}$ \ and \ $\bw \in \RR^d$.
\ Hence, by the proof of Theorem \ref{moment_m}, for each \ $t \in \RR_+$,
 \ $k \in \{1, \ldots, q\}$ \ and \ $\bw \in \RR^d$, \ there exists a
 polynomial \ $Q_{t,k,\bw} : \RR^d \to \RR$ \ having degree at most \ $k$ \ such
 that
 \[
   \EE\left[ \langle \bw, \bX_t \rangle^k \right]
   = \EE\bigl[Q_{t,j,\bw}(\bX_0)\bigr] ,
 \]
 where the coefficients of the polynomial \ $Q_{t,k,\bw}$ \ depends on
 \ $d$, $\bc$, $\Bbeta$, $\bB$, $\nu$, $\mu_1$, \ldots, $\mu_d$.

For all \ $a_1 , \ldots, a_k \in \RR$, \ we have
 \[
   a_1 \cdots a_k
   = \frac{1}{k! 2^k}
     \sum_{\ell_1=0}^1 \ldots \sum_{\ell_k=0}^1
      (-1)^{\ell_1+\cdots+\ell_k}
      \left[ (-1)^{\ell_1} a_1 + \cdots + (-1)^{\ell_k} a_k \right]^k .
 \]
Indeed, applying the multinomial theorem,
 \begin{gather*}
  \sum_{\ell_1=0}^1 \ldots \sum_{\ell_k=0}^1
   (-1)^{\ell_1+\cdots+\ell_k}
   \left[ (-1)^{\ell_1} a_1 + \cdots + (-1)^{\ell_k} a_k \right]^k \\
  = \sum_{\ell_1=0}^1 \ldots \sum_{\ell_k=0}^1
      (-1)^{\ell_1+\cdots+\ell_k}
      \sum_{\underset{j_1,\ldots,j_k \in \ZZ_+}{j_1+\cdots+j_k=k,}}
       \frac{k!}{j_1!\cdots j_k!}
       ((-1)^{\ell_1} a_1)^{j_1} \cdots ((-1)^{\ell_k} a_k)^{j_k}
   = S_1 + S_2 ,
 \end{gather*}
 where
 \begin{align*}
  S_1 &:= \sum_{\ell_1=0}^1 \ldots \sum_{\ell_k=0}^1
          (-1)^{\ell_1+\cdots+\ell_k}
          k! (-1)^{\ell_1} a_1 \cdots (-1)^{\ell_k} a_k , \\
  S_2 &:= \sum_{\ell_1=0}^1 \ldots \sum_{\ell_k=0}^1
          (-1)^{\ell_1+\cdots+\ell_k}
          \sum_{\underset{j_1,\ldots,j_k \in \ZZ_+}{j_1+\cdots+j_k=k,\,j_1\cdots j_k=0}}
           \frac{k!}{j_1!\cdots j_k!}
           ((-1)^{\ell_1} a_1)^{j_1} \cdots ((-1)^{\ell_k} a_k)^{j_k} .
 \end{align*}
Clearly \ $S_1 = 2^k k! a_1 \cdots a_k$, \ and \ $S_2 = 0$ \ because of
 cancellations.
Hence
 \begin{align*}
  &\EE\left( X_{t,i_1} \cdots X_{t,i_k} \right) \\
  &= \frac{1}{k! 2^k}
     \sum_{\ell_1=0}^1 \ldots \sum_{\ell_k=0}^1
      (-1)^{\ell_1+\cdots+\ell_k}
      \EE\left[\langle (-1)^{\ell_1} \be_{i_1} + \cdots + (-1)^{\ell_k} \be_{i_k}, \bX_t \rangle^k\right] \\
  &= \frac{1}{k! 2^k}
     \sum_{\ell_1=0}^1 \ldots \sum_{\ell_k=0}^1
      (-1)^{\ell_1+\cdots+\ell_k}
      \EE\bigl[Q_{t,k,(-1)^{\ell_1} \be_{i_1} + \cdots + (-1)^{\ell_k} \be_{i_k}}(\bX_0)
         \bigr]
  =: \EE[Q_{t,k,i_1,\ldots,i_k}(\bX_0)] ,
 \end{align*}
 which implies the statement.
\proofend

For central moments, we have the following recursion.

\begin{Thm}\label{central_moments_m}
Let \ $(\bX_t)_{t\in\RR_+}$ \ be a CBI process with parameters
 \ $(d, \bc, \Bbeta, \bB, \nu, \bmu)$ \ such that
 \ $\EE(\|\bX_0\|^q) < \infty$ \ and the moment conditions
 \eqref{moment_condition_m} hold with some \ $q \in \NN$.
\ Then
 \begin{equation}\label{help9}
  \begin{aligned}
   &\EE\left[(X_{t,j} - \EE(X_{t,j}))^k\right] \\
   &= k (k - 1)
      \sum_{i=1}^d
       c_i
       \int_0^t
        (\be_j^\top \ee^{(t-s)\tbB} \be_i)^2
        \EE\bigl[(\be_j^\top \ee^{(t-s)\tbB} (\bX_s - \EE(\bX_s)))^{k-2}
                 X_{s,i}\bigr]
        \dd s \\
   &\quad
    + \sum_{\ell=0}^{k-2}
         \binom{k}{\ell}
         \sum_{i=1}^d
          \int_0^t \int_{U_d}
           (\be_j^\top \ee^{(t-s)\tbB} \bz)^{k-\ell}
           \EE\Bigl[\bigl(\be_j^\top \ee^{(t-s)\tbB}(\bX_s - \EE(\bX_s))\bigr)^\ell
                    X_{s,i}\Bigr]
           \dd s \, \mu_i(\dd\bz) \\
   &\quad
    + \sum_{\ell=0}^{k-2}
         \binom{k}{\ell}
         \int_0^t \int_{U_d}
          (\be_j^\top \ee^{(t-s)\tbB} \bz)^{k-\ell}
          \EE\Bigl[\bigl(\be_j^\top \ee^{(t-s)\tbB}
                         (\bX_s - \EE(\bX_s))\bigr)^\ell\Bigr]
          \dd s \, \nu(\dd\bz)
  \end{aligned}
 \end{equation}
 for all \ $k \in \{1, \ldots, q\}$, \ $j \in \{1, \ldots, d\}$ \ and
 \ $t \in \RR_+$.
\ Moreover, for each \ $t \in \RR_+$, \ $k \in \{1, \ldots, q\}$ \ and
 \ $j \in \{1, \ldots, d\}$, \ there exists a polynomial
 \ $P_{t,k,j} : \RR^d \to \RR$ \ having degree at most \ $\lfloor k/2 \rfloor$
 \ such that
 \begin{align}\label{help9_polinom}
   \EE\left[ (X_{t,j} - \EE(X_{t,j}))^k \right]
   = \EE\left[ P_{t,k,j}(\bX_0) \right] , \qquad t \in \RR_+ .
 \end{align}
The coefficients of the polynomial \ $P_{t,k,j}$ \ depend on
 \ $d$, $\bc$, $\Bbeta$, $\bB$, $\nu$, $\mu_1$, \ldots, $\mu_d$.
\end{Thm}

\begin{Rem}
Note that in case of \ $\EE(\bX_t) = \bzero$, \ $t \in \RR_+$, \ formulae
 \eqref{SDE_truncated3} and \eqref{help9} coincide.
Indeed, if \ $\EE(\bX_t) = \bzero$, \ $t \in \RR_+$, \ then, by \eqref{EXbX},
 we have
 \[
   \EE(\be_j^\top\ee^{t\tbB} \bX_0)
   + \int_0^t \be_j^\top\ee^{u\tbB} \tBbeta \, \dd u
   = 0 , \qquad t \in \RR_+ , \qquad j \in \{1, \ldots, d\} .
 \]
Since \ $\be_j^\top\ee^{t\tbB} \bX_0$ \ is a non-negative random variable and
 \ $\RR_+ \ni t \mapsto \be_j^\top \ee^{t\tbB} \tBbeta$ \ is a non-negative
 continuous function, we obtain \ $\PP(\be_j^\top \ee^{t\tbB} \bX_0 = 0) = 1$
 \ and \ $\be_j^\top \ee^{t\tbB} \tBbeta = 0$ \ for all \ $t \in \RR_+$.
\ Consequently,
 \begin{align*}
  \EE\left[(\be_j^\top \ee^{t\tbB} \bX_0)^k\right]
      + k \int_0^t
           (\be_j^\top \ee^{(t-s)\tbB} \tBbeta)
           \EE\left[(\be_j^\top \ee^{(t-s)\tbB} \bX_s)^{k-1}\right]
           \dd s = 0
 \end{align*}
 for all \ $t \in \RR_+$, \ $j \in \{1, \ldots, d\}$ \ and \ $k \in \NN$,
 \ which yields that formulae \eqref{SDE_truncated3} and \eqref{help9}
 coincide.
\proofend
\end{Rem}

\noindent
\textbf{Proof of Theorem \ref{central_moments_m}.} \
Consider objects \textup{(E1)--(E4)} with initial value
 \ $\bxi = \by = (y_1, \ldots, y_d)^\top \in \RR_+^d$.
\ For each \ $K \in \NN$, \ let \ $(\bY_{K,t})_{t \in \RR_+}$ \ be a pathwise
 unique \ $\RR_+^d$-valued strong solution to the SDE \eqref{SDE_X_K} with
 initial value \ $\by$.
\ Using \eqref{Ito}, we obtain
 \begin{align*}
  \bw^\top \ee^{-t\tbB_K} (\bY_{K,t} - \EE(\bY_{K,t}))
  &= \sum_{i=1}^d
      \int_0^t
       \bw^\top \ee^{ -s\tbB_K} \be_i \sqrt{2 c_i Y_{K,s,i}} \, \dd W_{s,i} \\
  &\quad
     + \int_0^t \int_V
        \bw^\top \ee^{-s\tbB_K} h(\bY_{K,s-}, \br) \, \tN_K(\dd s, \dd\br)
 \end{align*}
 for all \ $\bw \in \RR^d$ \ and \ $t \in \RR_+$.
\ By the method of the proof of Theorem \ref{moment_m}, for a CBI process
 \ $(\bY_t)_{t\in\RR_+}$ \ having parameters
 \ $(d, \bc, \Bbeta, \bB, \nu, \bmu)$ \ with initial value \ $\by$, \ one can
 derive
 \begin{align}\label{help_central}
  \begin{split}
  &\EE\left[(\bw^\top \ee^{-t\tbB} (\bY_t - \EE(\bY_t))^k\right] \\
  &\quad = k (k - 1)
      \sum_{i=1}^d
       c_i \EE\left( \int_0^t
                      (\bw^\top \ee^{-s\tbB} \be_i)^2
                      (\bw^\top \ee^{-s\tbB} (\bY_s - \EE(\bY_s))^{k-2}
                      \, Y_{s,i}
                      \, \dd s \right) \\
  &\quad \phantom{=}
    + \sum_{\ell=0}^{k-2}
       \binom{k}{\ell}
       \EE\left( \int_0^t \int_V
                  (\bw^\top \ee^{-s\tbB} (\bY_s - \EE(\bY_s))^\ell
                  \bigl(\bw^\top \ee^{-s\tbB} h(\bY_s, \br)\bigr)^{k-\ell}
                  \, \dd s \, m(\dd\br) \right)
  \end{split}
 \end{align}
 for all \ $k \in \{2, \ldots, q\}$, \ where
 \begin{align*}
  &\EE\left( \int_0^t \int_V
                  (\bw^\top \ee^{-s\tbB} (\bY_s - \EE(\bY_s))^\ell
                  \bigl(\bw^\top \ee^{-s\tbB} h(\bY_s, \br)\bigr)^{k-\ell}
                  \, \dd s \, m(\dd\br) \right) \\
  &= \EE\left( \int_0^t \int_{U_d}
                (\bw^\top \ee^{-s\tbB} (\bY_s - \EE(\bY_s))^\ell
                \bigl(\bw^\top \ee^{-s\tbB} \br\bigr)^{k-\ell}
                \, \dd s \, \nu(\dd\br) \right) \\
  &\quad
     + \sum_{i=1}^d
        \EE\left( \int_0^t \int_{U_d} \int_{U_1}
                   (\bw^\top \ee^{-s\tbB} (\bY_s - \EE(\bY_s))^\ell
                   \bigl(\bw^\top \ee^{-s\tbB} \bz
                         \bbone_{\{s\leq Y_{s,i}\}}\bigr)^{k-\ell}
                   \, \dd s \, \mu_i(\dd\bz) \, \dd u \right)
 \end{align*}
 \begin{align*}
  &= \int_0^t \int_{U_d}
      \bigl(\bw^\top \ee^{-s\tbB} \br\bigr)^{k-\ell}
      \EE\bigl[(\bw^\top \ee^{-s\tbB} (\bY_s - \EE(\bY_s))^\ell \bigr]
      \, \dd s \, \nu(\dd\br) \\
  &\quad
     + \sum_{i=1}^d
        \int_0^t \int_{U_d}
         \bigl(\bw^\top \ee^{-s\tbB} \bz\bigr)^{k-\ell}
         \EE\bigl[(\bw^\top \ee^{-s\tbB} (\bY_s - \EE(\bY_s))^\ell Y_{s,i}\bigr]
         \, \dd s \, \mu_i(\dd\bz) .
 \end{align*}
As in the proof of Theorem \ref{moment_m}, this yields that the recursion
 \eqref{help9} holds for \ $\bY$, \ and, by the law of total probability, we
 obtain \eqref{help9} for \ $\bX$ \ as well.

Now we turn to prove \eqref{help9_polinom}.
As it was explained before, by the law of total probability, it is enough to
 prove \eqref{help9_polinom} for \ $\bY$.
\ Using the recursion \eqref{help9}, we obtain the existence of suitable
 polynomials \ $P_{t,k,j}$, \ $t \in \RR_+$ \ $k \in \{1, \ldots, q\}$,
 \ $j \in \{1, \ldots, d\}$, \ by induction with respect to \ $k$.
\ Indeed, for \ $k = 1$, \ we have \ $\EE[Y_{t,j} - \EE(Y_{t,j})] = 0$,
 \ $j \in \{1, \ldots, d\}$, \ $t \in \RR_+$.
\ For \ $k = 2$, \ by \eqref{help9}, we have
 \begin{align}\label{central_moment_2}
  \begin{split}
   &\EE\left[(Y_{t,j} - \EE(Y_{t,j}))^2\right]
    = 2 \sum_{i=1}^d
        c_i \int_0^t (\be_j^\top \ee^{(t-s)\tbB} \be_i)^2  \EE(Y_{s,i}) \, \dd s \\
   &\hspace*{14.7mm}
    + \sum_{i=1}^d
       \int_0^t \int_{U_d}
        (\be_j^\top \ee^{(t-s)\tbB} \bz)^2 \EE(Y_{s,i})
        \, \dd s \, \mu_i(\dd\bz)
    + \int_0^t \int_{U_d}
       (\be_j^\top \ee^{(t-s)\tbB} \bz)^2 \dd s \, \nu(\dd\bz)
  \end{split}
 \end{align}
 for all \ $j \in \{1, \ldots, d\}$ \ and \ $t \in \RR_+$.
\ Thus \ $\EE\left[(Y_{t,j} - \EE(Y_{t,j}))^2\right] = P_{t,2,j}(\by)$,
 \ where \ $P_{t,2,j} : \RR^d \to \RR$ \ is a polynomial of degree at most 1,
 since
 \ $\EE(Y_{s,i})
    = \be_i^\top \ee^{s\tbB} \by
      + \int_0^s \be_i^\top \ee^{u\tbB} \tBbeta \, \dd u$,
 \ $s \in \RR_+$, \ from \eqref{EXbX}.

Now, suppose that for some \ $k' \in \NN$ \ with \ $2 k' + 1 \leq q$,
 \ suitable polynomials \ $P_{t,1,j}$, \ldots, $P_{t,2k',j}$ \ exist for all
 \ $t \in \RR_+$ \ and \ $j \in \{1, \ldots, d\}$.
\ We apply the recursion \eqref{help9} for \ $k = 2 k' + 1$.
\ Then for each \ $\ell \in \{0, 1, \ldots, 2k' - 1\}$ \ and
 \ $s, t \in \RR_+$ \ with \ $s \leq t$, \ the function
 \[
   \RR_+^d \ni \by
   \mapsto \EE\Bigl[ \bigl(\be_j^\top \ee^{(t-s)\tbB}
                          (\bY_s - \EE(\bY_s))\bigr)^\ell \Bigr]
 \]
 is a polynomial of degree at most
 \ $\lfloor \ell/2 \rfloor \leq \lfloor (2k' - 1)/2 \rfloor = k' - 1$.
\ Moreover, for each \ $\ell \in \{0, 1, \ldots, 2 k' - 1\}$ \ and
 \ $s, t \in \RR_+$ \ with \ $s \leq t$, \ the function
 \[
   \RR_+^d \ni \by
   \mapsto \EE\Bigl[ \bigl(\be_j^\top \ee^{(t-s)\tbB}
                          (\bY_s - \EE(\bY_s))\bigr)^\ell Y_{s,i} \Bigr]
 \]
 is a polynomial of degree at most
 \ $\max\{\lfloor \ell/2 \rfloor +1, \lfloor (\ell + 1)/2 \rfloor \}
    \leq \max\{ k', \lfloor (2 k')/2 \rfloor\} = k'$,
 \ since, by \eqref{EXbX},
 \[
   Y_{s,j} = \EE(Y_{s,j}) + (Y_{s,j} - \EE(Y_{s,j}))
          = \be_j^\top \ee^{s\tbB} \by
            + \int_0^s \be_j^\top \ee^{v\tbB} \tBbeta \, \dd v
            + (Y_{s,j} - \EE(Y_{s,j})) .
 \]
Consequently, by \eqref{help9},
 \ $\RR_+^d \ni \by \mapsto \EE\left[(Y_{t,j} - \EE(Y_{t,j}))^{2k'+1} \right]$
 \ is a polynomial of degree at most \ $k' = \lfloor (2k' + 1)/2 \rfloor$,
 \ and we conclude the existence of suitable polynomials \ $P_{t,2k'+1,j}$ \ for
 all \ $t \in \RR_+$ \ and \ $j \in \{1, \ldots, d\}$.

In a similar way, if for some \ $k' \in \NN$ \ with \ $2 k' + 2 \leq q$,
 \ suitable polynomials \ $P_{t,1,j}$, \ldots, $P_{t,2k'+1,j}$ \ exist for all
 \ $t \in \RR_+$ \ and \ $j \in \{1, \ldots, d\}$, \ then we apply the
 recursion \eqref{help9} for \ $k = 2 k' + 2$.
\ Then for each \ $\ell \in \{0, 1, \ldots, 2k'\}$, \ the function
 \ $\RR_+^d \ni \by
    \mapsto \EE\Bigl[ \bigl(\be_j^\top \ee^{(t-s)\tbB}
                           (\bY_s - \EE(\bY_s))\bigr)^\ell \Bigr]$
 \ is a polynomial of degree at most
 \ $\lfloor \ell/2 \rfloor \leq \lfloor (2k')/2 \rfloor = k'$.
\ Further, for each \ $\ell \in \{0, 1, \ldots, 2k'\}$, \ the function
 \ $\RR_+^d \ni \by \mapsto
    \EE\Bigl[ \bigl(\be_j^\top \ee^{(t-s)\tbB}
                   (\bY_s - \EE(\bY_s))\bigr)^\ell Y_{s,i} \Bigr]$
 \ is a polynomial of degree at most
 \ $\max\{ \lfloor \ell/2 \rfloor + 1, \lfloor (\ell + 1)/2 \rfloor \}
    \leq \max\{ k'+1, \lfloor (2k' + 1)/2 \rfloor \}= k' + 1$.
\ Consequently, by \eqref{help9},
 \ $\RR_+^d \ni \by \mapsto \EE\left[(Y_{t,j} - \EE(Y_{t,j}))^{2k'+2} \right]$
 \ is a polynomial of degree at most \ $k' + 1 = \lfloor (2k' + 2)/2 \rfloor$,
 \ and we conclude the existence of suitable polynomials \ $P_{t,2k'+2,j}$ \ for
 all \ $t \in \RR_+$ \ and \ $j \in \{1, \ldots, d\}$.
\proofend

For mixed central moments, we have the following corollary.

\begin{Cor}\label{Cor_mixed_central_moments}
Let \ $(\bX_t)_{t\in\RR_+}$ \ be a CBI process with parameters
 \ $(d, \bc, \Bbeta, \bB, \nu, \bmu)$ \ such that
 \ $\EE(\|\bX_0\|^q) < \infty$ \ and the moment conditions
 \eqref{moment_condition_m} hold with some \ $q \in \NN$.
\ Then for all \ $t \in \RR_+$, \ $k \in \{1, \ldots, q\}$ \ and
 \ $i_1, \ldots, i_k \in \{1, \ldots, d\}$, \ there exists
 a polynomial \ $P_{t,k,i_1,\ldots,i_k} : \RR^d \to \RR$ \ having degree at most
 \ $\lfloor k/2 \rfloor$ \ such that
 \begin{align}\label{polinomP}
  \EE\big[ (X_{t,i_1} - \EE(X_{t,i_1}) )\cdots (X_{t,i_k} - \EE(X_{t,i_k})) \big]
  = \EE(P_{t,k,i_1,\ldots,i_k}(\bX_0)).
 \end{align}
The coefficients of the polynomial \ $P_{t,k,i_1,\ldots,i_k}$ \ depend on
 \ $d$, $\bc$, $\Bbeta$, $\bB$, $\nu$, $\mu_1$, \ldots, $\mu_d$.
\end{Cor}

\noindent{\bf Proof.}
Replacing \ $\bw$ \ by \ $\ee^{t\tbB^\top} \!\! \bw$ \ in \eqref{help_central},
 and then using the law of total probability, one obtains
 \begin{equation*}
  \begin{aligned}
   &\EE\left[\langle \bw, \bX_t - \EE(\bX_t) \rangle^k\right] \\
   &= k (k - 1)
      \sum_{i=1}^d
       c_i
       \int_0^t
        (\bw^\top \ee^{(t-s)\tbB} \be_i)^2
        \EE\bigl[(\bw^\top \ee^{(t-s)\tbB} (\bX_s - \EE(\bX_s)))^{k-2}
                 X_{s,i}\bigr]
        \dd s \\
   &\quad
    + \sum_{\ell=0}^{k-2}
         \binom{k}{\ell}
         \sum_{i=1}^d
          \int_0^t \int_{\cU_d}
           (\bw^\top \ee^{(t-s)\tbB} \bz)^{k-\ell}
           \EE\Bigl[\bigl(\bw^\top \ee^{(t-s)\tbB}(\bX_s - \EE(\bX_s))\bigr)^\ell
                    X_{s,i}\Bigr]
           \dd s \, \mu_i(\dd\bz) \\
   &\quad
    + \sum_{\ell=0}^{k-2}
         \binom{k}{\ell}
         \int_0^t \int_{\cU_d}
          (\bw^\top \ee^{(t-s)\tbB} \bz)^{k-\ell}
          \EE\Bigl[\bigl(\bw^\top \ee^{(t-s)\tbB}
                         (\bX_s - \EE(\bX_s))\bigr)^\ell\Bigr]
          \dd s \, \nu(\dd\bz)
  \end{aligned}
 \end{equation*}
 for all \ $t \in \RR_+$, \ $k \in \{1, \ldots, q\}$ \ and \ $\bw \in \RR^d$,
 \ and hence, by the proof of Theorem \ref{central_moments_m},
  for each \ $t \in \RR_+$, \ $k \in \{1, \ldots, q\}$ \ and
 \ $\bw \in \RR^d$, \ there exists a polynomial \ $P_{t,k,\bw} : \RR^d \to \RR$
 \ having degree at most \ $\lfloor k/2 \rfloor$, \ such that
 \[
   \EE\left[ \langle \bw, \bX_t - \EE(\bX_t) \rangle^k \right]
   = \EE\bigl[P_{t,k,\bw}(\bX_0)\bigr] ,
 \]
 where the coefficients of the polynomial \ $P_{t,k,\bw}$ \ depend on
 \ $d$, $\bc$, $\Bbeta$, $\bB$, $\nu$, $\mu_1$, \ldots, $\mu_d$.
\ The proof can be finished as the proof of Corollary \ref{Cor_mixed_moments}.
\proofend

\begin{Pro}\label{moment_formula_2}
Let \ $(\bX_t)_{t\in\RR_+}$ \ be a CBI process with parameters
 \ $(d, \bc, \Bbeta, \bB, \nu, \bmu)$ \ such that \ $\EE(\|\bX_0\|^2) < \infty$
 \ and the moment conditions \eqref{moment_condition_m} hold with \ $q = 2$.
\ Then for all \ $t \in \RR_+$, \ we have
 \begin{align*}
  \var(\bX_t)
  &= \sum_{\ell=1}^d
      \int_0^t
       (\be_\ell^\top \ee^{(t-u)\tbB} \EE(\bX_0))
       \ee^{u\tbB} \bC_\ell \ee^{u\tbB^\top}
       \dd u
     + \int_0^t
        \ee^{u\tbB}
        \left( \int_{U_d} \bz \bz^\top \nu(\dd \bz) \right)
        \ee^{u\tbB^\top}
        \dd u \\
  &\quad
     + \sum_{\ell=1}^d
        \int_0^t
         \left( \int_0^{t-u} \be_ \ell^\top \ee^{v\tbB} \tBbeta \, \dd v \right)
         \ee^{u\tbB} \bC_\ell \ee^{u\tbB^\top}
         \dd u ,
 \end{align*}
 where
 \[
   \bC_\ell := 2 c_\ell \be_\ell \be_\ell^\top
              + \int_{U_d} \bz \bz^\top \mu_\ell(\dd \bz)
           \in \RR_+^{d \times d} , \qquad
   \ell \in \{1, \ldots, d\} .
 \]
\end{Pro}

\noindent
\textbf{Proof.} \
By \eqref{central_moment_2}, we have
 \begin{multline*}
  \be_j^\top \EE\left[(\bX_t - \EE(\bX_t)) (\bX_t - \EE(\bX_t))^\top\right] \be_j
  = \be_j^\top \var(\bX_t) \be_j
  = \EE\left[(X_{t,j} - \EE(X_{t,j}))^2\right] \\
  = \sum_{\ell=1}^d
     \int_0^t
      \be_j^\top \ee^{(t-u)\tbB} \bC_\ell \, \ee^{(t-u)\tbB^\top} \be_j
      \EE(X_{u,\ell})
      \, \dd u
    + \int_0^t
       \be_j^\top \ee^{(t-u)\tbB}
       \left( \int_{U_d} \bz \bz^\top \, \nu(\dd\bz) \right)
       \ee^{(t-u)\tbB^\top} \be_j \, \dd u ,
 \end{multline*}
 which is finite by \eqref{moment_condition_m} with \ $q = 2$ \ and part (v)
 of Definition \ref{Def_admissible}.
Using the identities
 \[
   \be_i^\top \var(\bX_t) \be_j
   = \frac{1}{4}
     \left[ (\be_i + \be_j)^\top \var(\bX_t) (\be_i + \be_j)
            - (\be_i - \be_j)^\top \var(\bX_t) (\be_i - \be_j) \right]
 \]
 for \ $i,j\in\{1,\ldots,d\}$, \ and
 \ $\var(\bX_t)
    = \sum_{i=1}^d \sum_{j=1}^d \be_i (\be_i^\top \var(\bX_t) \be_j) \be_j^\top$,
 \ we obtain
 \[
   \var(\bX_t)
   = \sum_{\ell=1}^d
      \int_0^t
       \ee^{(t-u)\tbB} \bC_\ell \, \ee^{(t-u)\tbB^\top} \EE(X_{u,\ell}) \, \dd u
    + \int_0^t
       \ee^{(t-u)\tbB}
       \left( \int_{U_d} \bz \bz^\top \, \nu(\dd\bz) \right)
       \ee^{(t-u)\tbB^\top} \dd u .
 \]
By \eqref{EXbX}, we have
 \ $\EE(X_{u,\ell})
    = \be_\ell^\top \ee^{u\tbB} \EE(\bX_0)
      + \int_0^u \be_\ell^\top \ee^{v\tbB} \tBbeta \, \dd v$,
 \ thus
 \begin{align*}
  \var(\bX_t)
  &= \sum_{\ell=1}^d \int_0^t
      (\be_\ell^\top \ee^{u\tbB} \EE(\bX_0))
      \, \ee^{(t-u)\tbB} \bC_\ell \ee^{(t-u)\tbB^\top}
      \dd u
     + \int_0^t
        \ee^{u\tbB} \left( \int_{U_d} \bz \bz^\top \nu(\dd \bz) \right)
        \ee^{u\tbB^\top} \dd u \\
  &\phantom{=\;}
     + \sum_{\ell=1}^d
        \int_0^t
         \left(\int_0^u \be_\ell^\top \ee^{v\tbB}\tBbeta \,\dd v
         \ee^{(t-u)\tbB} \bC_\ell \ee^{(t-u)\tbB^\top} \right) \dd u \\
  &= \sum_{\ell=1}^d
      \int_0^t
       (\be_\ell^\top \ee^{(t-v)\tbB} \EE(\bX_0))
       \, \ee^{v\tbB} \bC_\ell \ee^{v\tbB^\top}
       \, \dd v
     + \int_0^t
        \ee^{u\tbB} \left( \int_{U_d} \bz \bz^\top \nu(\dd \bz) \right)
        \ee^{u\tbB^\top} \dd u \\
  &\phantom{=\;}
     + \sum_{\ell=1}^d
        \int_0^t
         \left(\int_0^{t-u} \be_\ell^\top \ee^{v\tbB} \tBbeta \, \dd v\right)
         \ee^{u\tbB} \bC_\ell \ee^{u\tbB^\top} \dd u ,
 \end{align*}
 and hence we obtain the statement.
\proofend

\section*{Acknowledgements}
We would like to thank the referee for his/her comments that helped us to improve the presentation
 of the paper.

\end{document}